\documentclass[letterpaper,11pt]{article}

\usepackage{amssymb}

\usepackage{amsthm}
\usepackage{amsmath}
\usepackage[utf8x]{inputenc}
\usepackage[T1]{fontenc}

\usepackage[hidelinks,hypertexnames=false]{hyperref}
\usepackage[numbers]{natbib}
\usepackage[margin=1in]{geometry}

\usepackage{chngcntr}
\counterwithin{table}{subsection}

\usepackage{graphicx} 
\usepackage{epstopdf}  \usepackage{enumerate}
\usepackage{url} \usepackage{verbatim} \usepackage{float} \usepackage{subfig}

\usepackage{tabularx,ragged2e,booktabs,caption} \usepackage{makecell}
\usepackage[boxruled]{algorithm2e}

\usepackage{comment}
\usepackage{enumitem}

\title{Partial facial reduction: simplified, equivalent SDPs \\ 
via  approximations of the PSD cone}

\DeclareMathOperator*{\bigO}{\ensuremath{\mathcal{O}}}

\DeclareMathOperator*{\rank}{rank}
\DeclareMathOperator*{\cprank}{rank_{cp}}

\DeclareMathOperator*{\range}{range}
\DeclareMathOperator*{\nullspace}{null}
\DeclareMathOperator*{\trace}{Tr}

\DeclareMathOperator*{\lin}{lin}

\DeclareMathOperator*{\relint}{relint}
\DeclareMathOperator*{\nnz}{nnz}
\DeclareMathOperator*{\vect}{vect}

\newcommand{\setName}{\mathbb{W}}
\newcommand{\coneName}{\mathcal{K}}
\newcommand{\subspaceName}{\mathcal{A}}
\newcommand{\faceName}{\mathcal{F}}
\newcommand{\approxName}{\mathcal{C}(\setName)^{*}}
\newcommand{\outerName}{\mathcal{C}(\setName)}

\newcommand{\timeLP}{t_{LPs}}

\newcommand{\FW}{\mathcal{FW}}
\newcommand{\SDD}{\mathcal{SDD}}
\newcommand{\DD}{\mathcal{DD}}
\newcommand{\D}{\mathcal{D}}

\newcommand{\perpName}{\mathcal{M}}

\author{Frank Permenter \and Pablo Parrilo \thanks{ The authors are with the Laboratory for Information and Decision Systems, Department of Electrical Engineering and Computer Science, Massachusetts Institute of Technology, Cambridge MA 02139, USA. Email:
\{fperment,parrilo\}@mit.edu.}%
}

\begin{document}
\newtheorem{thm}{Theorem}
\newtheorem{cor}{Corollary}
\newtheorem{lem}{Lemma}
\newtheorem{fact}{Fact}
\newtheorem{rem}{Remark}
\newtheorem{prob}{Problem}
\newtheorem{prop}{Proposition}
\newtheorem{defn}{Definition}
\newtheorem{ex}{Example}
\newtheorem{cond}{Condition}

\maketitle

\begin{abstract}

We develop a practical semidefinite programming (SDP) facial reduction procedure that utilizes computationally efficient  approximations of the positive semidefinite cone. The proposed method simplifies SDPs with no strictly feasible solution (a frequent output of parsers) by solving a sequence of easier optimization problems and could be a useful pre-processing technique for SDP solvers.  We demonstrate  effectiveness of the method on SDPs arising in practice, and describe our publicly-available software implementation. We also show how to find maximum rank matrices in our PSD cone approximations (which helps us find maximal simplifications), and we give a post-processing procedure for dual solution recovery that generally applies to facial-reduction-based  pre-processing techniques.  Finally, we show how   approximations can be chosen to preserve problem sparsity.

\end{abstract}

\section{Introduction}
The feasible set of a semidefinite program (SDP) is described by the intersection of an affine subspace with the   cone of matrices that are positive semidefinite (PSD).
In practice, this intersection may contain no matrices that are strictly positive definite, i.e., strict feasibility may fail.  This is problematic for two reasons.  One,
strong duality is not guaranteed. Two, the SDP (if feasible) is unnecessarily large in the sense it can be reformulated using a
smaller PSD cone and a lower dimensional subspace. To see this latter point, consider
the following motivating example: 
\paragraph{Motivating example}
\begin{align*}
\begin{array}{l}
\begin{array}{ll}
\mbox{Find} &  y_1,y_2,y_3 \in \mathbb{R}\\
 \mbox{{subject to  }} & 
\end{array} \\
\;\; \;\; \subspaceName(y) =  \left(\begin{array}{ccc} y_1 & 0 & 0 \\ 
0 & - y_1 & y_2 \\
0 & y_2 & y_2 + y_3 \\ 
\end{array}\right) \succeq 0.
\end{array} 
\end{align*}
Taking $v = (1,1,0)^T$, it is clear that $v^T \subspaceName(y) v = 0$ independent of $(y_1,y_2,y_3)$.  In other words, there is no $(y_1,y_2,y_3)$ for which $\subspaceName(y)$ is positive definite.
It also holds that $(y_1,y_2,y_3)$ is a feasible point of the above SDP if and only if it is a feasible point of
\begin{align*}
\begin{array}{ll}
\mbox{Find} &  y_1,y_2,y_3 \in \mathbb{R}\\
 \mbox{{subject to  }} &  y_1 = y_2 = 0, y_3 \ge 0.
\end{array} 
\end{align*}
In other words,  the above $3 \times 3$ semidefinite constraint is equivalent to linear equations and a linear inequality (i.e. a $1 \times 1$ semidefinite constraint).

While SDPs of this type may seem rare in practice, they are a frequent output of \emph{parsers} (e.g. \cite{YALMIP}, \cite{sostools}) used, for example, to formulate
 SDP-based relaxations of algebraic problems. In some cases, these SDPs arise because   the parser does not exploit available problem structure  
(cf. the SDP relaxations of graph partitioning problem \cite{wolkowicz1999semidefinite}, where problem structure is carefully exploited to ensure strict feasibility). In other situations, 
all relevant structure is not apparent from the problem's natural high level description  (which motivates the post-processing of solutions in \cite{lofberg2009pre}).
Thus, based on their prevalence,  checking for and simplifying an SDP of this type is    a  practically useful \emph{pre-processing} step, assuming it can be done efficiently.

To check for (and simplify) such an SDP, one can  execute the  facial reduction algorithm of \citet{borwein1981regularizing}, or the simplified versions of \citet{pataki2013simple} and \citet{waki2013facial}, to 
 find a \emph{face} of the PSD cone  containing the feasible set. The desired simplifications are then obtained by  reformulating the SDP as an optimization problem over this face.
 Unfortunately, the problem of finding  a face is itself an SDP, which  may be too expensive to solve in the context of pre-processing.  In addition, reformulating the SDP accurately, in a way that preserves sparsity, can be difficult.

To address these issues, this paper presents a facial reduction algorithm modified in a simple way: rather than search over all possible faces, 
our method searches over just a subset defined by a specified \emph{approximation} of the PSD cone. As we show, the specified approximation allows one to control pre-processing effort, preserve sparsity, and  accurately reformulate a given SDP. Natural choices for the approximation are also effective in practice, as we illustrate with examples.

This paper is organized as follows.  In Section~\ref{sec:background}, we give a  short derivation of a facial reduction algorithm for general conic optimization problems 
using basic tools from convex analysis.
We then specialize this algorithm to semidefinite programs.  In Section~\ref{sec:approach}, we modify this  algorithm to yield our technique and describe example approximations of the PSD cone. We then show
how to find maximum rank solutions to conic optimization problems formulated over these approximations (which helps us find faces of minimal dimension).  Section~\ref{sec:simpleEx} shows
how to reformulate a given SDP over an identified face and illustrates how the chosen approximation affects sparsity.  Simple illustrative
examples   are then given.  Section~\ref{sec:dualRecov} discusses the   issue of dual solution recovery
and generalizes a recovery procedure described in \cite{pataki2016bad}.
 The results of this section are not specific to our modified facial reduction procedure
and are   relevant to other pre-processing techniques based on facial reduction. Section~\ref{sec:implementation} describes a freely-available implementation of our procedure
and Section~\ref{sec:examples} illustrates effectiveness of the method on examples arising in practice.

\subsection{Prior work}
 General algorithms for facial reduction include the original of \citet{borwein1981regularizing} and the simplified versions of \citet{pataki2013simple} and  \citet{waki2013facial}. Application specific approaches have also been developed: e.g., for SDPs arising in Euclidean distance matrix completion (\citet{krislock2010explicit}), protein structure identification 
(\citet{Alipanahi2012};  \citet{burkowski2011efficient}), graph partitioning 
 (\citet{wolkowicz1999semidefinite}),  quadratic assignment (\citet{zhao1998semidefinite}) and max-cut (\citet{anjos2002strengthened}). 
 \citet{waki2010facial} and the current authors \cite{permenter2014basis} also apply facial reduction to the problem of basis selection in sums-of-squares optimization.

In addition, facial reduction can be used to ensure that strong duality holds.  Indeed, this was the original motivation of  the  Borwein~and~Wolkowicz  algorithm \cite{borwein1981regularizing},
which given a feasible  optimization problem outputs one that  satisfies Slater's condition.
Facial reduction is also the  basis for so-called extended duals, which are generalized dual programs for which strong duality always holds.
Extended duals are studied by  \citeauthor{pataki2013simple} for optimization problems over \emph{nice} cones  \cite{pataki2013simple} and include   Ramana's dual for SDP \cite{ramana1997exact} \cite{ramana1997strong}.

A dual view of  facial reduction  was given  by Luo, Sturm, and  Zhang  \cite{Luo97dualityresults}.
There, the authors describe a so-called conic expansion algorithm 
that grows the dual cone to include additional linear functionals non-negative
on  the feasible set. In \cite{waki2013facial},  Waki and   Muramatsu give a
facial reduction procedure and explicitly  relate it to  conic 
expansion.

The idea of using facial reduction as a pre-processing step for SDP was described in \cite{gruber1998presolving} by Gruber et al.  The
authors note the expense of identifying lower dimensional faces as well as issues of numerical reliability that may arise. 
In \cite{cheung2011preprocessing},   Cheung,  Schurr, and  Wolkowicz address the issue of numerical reliability, giving a facial
reduction algorithm that identifies a \emph{nearby} problem in a backwards stable manner.  


Finally, our technique is consistent with a philosophy of \citeauthor{andersen1995presolving} put forth in \cite{andersen1995presolving}.
There, the authors argue the best strategy for pre-processing LPs is to find \emph{simple} simplifications \emph{quickly}.  The method we present is consistent with this philosophy in  that the specified approximation defines the notion of ``simple'' and its search complexity
defines the notion of ``quick.''

\subsection{Contributions}

 \paragraph{Partial facial reduction} Our main contribution  is a   pre-processing technique for semidefinite programs, based
on facial reduction, that allows one to specify pre-processing effort, preserve problem sparsity, and ensure accuracy of problem reformulations. Given any SDP,
the technique searches for an equivalent reformulation over a lower dimensional face, where a
user-specified approximation of the PSD cone   controls the size of this search space.  In addition, natural choices for the approximation preserve sparsity, as we explain in Section~\ref{sec:simpleEx}. Finally,  if a polyhedral approximation is specified, our method solves only linear programs, which are accurately solved both in theory and in practice (and in exact arithmetic, if desired).
 \paragraph{Maximum rank solutions} 
 Related to finding a face of minimal dimension is finding a maximum rank matrix in a subspace intersected  with a specified approximation.
We show (Corollary~\ref{cor:maxRank}) how to  find such a matrix when the approximation  equals the Minkowski sum of faces of the PSD cone.  
Approximations of this type include diagonally-dominant \cite{barker1975cones}, scaled diagonally-dominant, and factor-width-$k$   \cite{boman2005factor} approximations.   
 
 \paragraph{Dual solution recovery} We  give and study a simple algorithm for dual solution recovery (Algorithm~\ref{alg:frPrim_DualRecov}) that generalizes an approach from \cite{pataki2016bad}. Dual solution recovery is a critical post-processing step for primal-dual solvers, where are often agnostic to which problem---primal or dual---is of actual interest to a user. For this reason, recovery has received much attention in linear programming \cite{andersen1995presolving}. Our  recovery procedure applies generally to conic optimization problems pre-processed using facial reduction techniques; in other words, it is
not specific to SDP and does not depend on the approximations we introduce.  Since pre-processing may remove duality gaps, dual solution recovery is \emph{not always possible}. Hence, we give conditions  
(Conditions~\ref{cond:sdpRecov}~and~\ref{cond:recEnsure})  characterizing success of the procedure for  SDPs---the class of conic optimization problem of primary interest. 

 \paragraph{Software implementation}  We have implemented our technique  in MATLAB in a set of scripts we call {\tt frlib}, 
 available  at \tt{www.mit.edu/\textasciitilde fperment}\rm.
If interfaced directly, the code takes as input SDPs in SeDuMi format  \cite{sturm1999using}.  It can also be interfaced via the 
parser  YALMIP \cite{YALMIP}.

\section{Background on facial reduction} \label{sec:background}
In this section, we  define our notation, collect basic facts and definitions, and  describe faces of the PSD cone.
We then review facial reduction, giving a simple, self-contained derivation of the simplified facial reduction algorithm of \citet{pataki2013simple}. We then specialize this algorithm to semidefinite programs.  

\subsection{Notation and preliminaries}

Let $\mathcal{E}$ denote a finite-dimensional vector space over $\mathbb{R}$ with inner product $\langle \cdot, \cdot \rangle$.  For a subset $\mathcal{S}$ of $\mathcal{E}$, let $\lin \mathcal{S} \subseteq \mathcal{E}$ denote 
the linear span of elements in $\mathcal{S}$ and let $\mathcal{S}^{\perp} \subseteq \mathcal{E}$ denote the orthogonal complement of $\lin \mathcal{S}$.  For $y \in \mathcal{E}$, let $\lin y$ denote $\lin \{ y \}$ and let $y^{\perp}$  denote $\{ y \}^{\perp}$.
A  convex cone $\coneName$ is a  convex subset of $\mathcal{E}$ (not necessarily full dimensional) that satisfies
\[
 x  \in \coneName \Rightarrow  \alpha x \in \coneName \;\;\; \forall \alpha \ge 0.
\]
The dual cone of $\coneName$, denoted  $\coneName^*$, is the set of linear functionals non-negative on $\coneName$:
\[
  \coneName^* := \left\{ y : \langle y,x\rangle \ge 0 \;\;\; \forall x \in \coneName\right\}.
\]
A  face $\faceName$ of a convex cone $\coneName$ is a convex subset that satisfies 
\[
 \frac{a+b}{2} \in \faceName \mbox{ and } a,b \in \coneName \Rightarrow a,b \in \faceName.
\]
A face is \emph{proper} if it is non-empty and not equal to $\coneName$. Faces of convex cones are also convex cones, and the relation ``is a face of'' is  transitive; if $\mathcal{F}_2$ is a face of $\mathcal{F}_1$ and $\mathcal{F}_3$ is a face of $\mathcal{F}_2$,
 then $\mathcal{F}_3$ is face of $\mathcal{F}_1$.
For any $s \in \coneName^*$, the set $\coneName \cap s^{\perp}$ is a face of $\coneName$.
Further, if $\coneName$ is closed, it holds that $(\coneName \cap s^{\perp})^* = \overline{\coneName^* + \lin s}$ (where we let $\overline{\mathcal{S}}$ denote
the closure of a set $\mathcal{S}$).

For discussions specific to semidefinite programming, we let $\mathbb{S}^n$ denote the vector space of $n\times n$ symmetric matrices and $\mathbb{S}^n_{+} \subseteq \mathbb{S}^n$
denote the convex cone of matrices that are positive semidefinite. We will use capital letters to denote elements of
$\mathbb{S}^n$ to emphasize that they are matrices.  For $A,B \in \mathbb{S}^n$, we let $A \cdot B$  denote the trace inner product $\trace AB$. 
Finally, we let  $A \succeq 0$ (resp. $A \succ 0$)  denote the condition that $A$ is positive semidefinite (resp.  positive definite).

\subsection{Faces of $\mathbb{S}^n_{+}$} 
A set is a  face of $\mathbb{S}^n_{+}$ if and only if it equals the set of all $n \times n$ PSD matrices with range contained in a given $d$-dimensional subspace \cite{barker1975cones} \cite{pataki2000geometry}. 
 Using this fact, one can describe a proper face $\faceName$  (and the dual cone $\faceName^*$) 
using an invertible matrix $(U,V) \in \mathbb{R}^{n \times n}$, where the range of $U \in \mathbb{R}^{n \times d}$ equals 
this subspace and  the range of $V\in \mathbb{R}^{n \times n-d}$ equals $(\range U)^{\perp}$. We collect  such descriptions in the following.

\begin{lem} \label{lem:facesofPSD}
A non-zero, proper face of $\mathbb{S}^n_{+}$ is a set $\faceName$ of the form
\begin{align}
 \faceName :=&  \left\{ (U,V) \left( \begin{array}{cc} W & 0 \\ 0 & 0 \end{array} \right) (U,V)^T : W \in \mathbb{S}^d_{+} \right\} \label{eq:faceGen} \\
 =& \left\{ X \in \mathbb{S}^n : U^T X U \succeq 0, U^T X V = 0, V^T X V = 0 \right\}, \label{eq:faceCnt}
\end{align}
where $(U,V) \in \mathbb{R}^{n \times n}$ is an invertible matrix satisfying $U^T V = 0$  (i.e. $\range V = (\range U)^{\perp}$). Moreover, the dual cone $\faceName^*$ satisfies
\begin{align}
 \faceName^* &= \left\{ (U,V) \left( \begin{array}{cc} W & Z \\ Z^T & R \end{array} \right) (U,V)^T :  W \in \mathbb{S}^d_{+}, Z \in \mathbb{R}^{d \times n-d}, R \in \mathbb{S}^{n-d} \right\} \label{eq:faceDualGen} \\
  &= \left\{ X \in \mathbb{S}^n : U^T X U \in \mathbb{S}^d_{+} \right\}. \label{eq:faceDualCnt} 
\end{align}
\end{lem}
\noindent Here, \eqref{eq:faceGen} and \eqref{eq:faceDualGen} represent a face $\faceName$ and its dual cone $\faceName^*$ in terms of \emph{generators} whereas   \eqref{eq:faceCnt} and \eqref{eq:faceDualCnt} represent
these sets in terms of \emph{constraints}.  Either representation can be preferred depending on context.  Based off of \eqref{eq:faceGen}, we will often refer to a face $\faceName$ using the notation $U \mathbb{S}^d_{+} U^T$.

\subsection{Facial reduction of conic optimization problems}

\paragraph{Conic optimization problems}
The feasible set of a conic optimization problem is described by the intersection of an affine subspace $\subspaceName$ with 
a convex cone $\coneName$, where both $\subspaceName$ and $\coneName$ are subsets of the inner product space $\mathcal{E}$.  If one defines the affine subspace $\subspaceName$ in terms of a linear map $A : \mathbb{R}^m \rightarrow \mathcal{E}$ 
and a point $c \in \mathcal{E}$, i.e.
\[
\subspaceName  := \{ c - Ay : y \in \mathbb{R}^m \},
\]
one can express a  conic optimization problem as follows:
\[
         \mbox{maximize } b^Ty   \mbox{ subject to }    c - Ay  \in \coneName,
\]
where $b \in \mathbb{R}^m$ defines a linear objective function. This conic optimization problem is \emph{feasible} if
 $\subspaceName \cap \coneName$ is non-empty and \emph{strictly
   feasible} if $\subspaceName \cap \relint \coneName$ is
 non-empty.

\paragraph{Reformulation over a face}
 If a conic optimization problem is feasible but not strictly feasible, 
it can be reformulated as an optimization problem over a lower dimensional face of $\coneName$.  This fact will follow from the following lemma, which holds for arbitrary convex cones.  See also Lemma 1 of \cite{pataki2013simple},   Theorem 7.1 of \cite{borwein1981regularizing}, Lemma 12.6 of \cite{cheung2011preprocessing}, and Lemma 3.2 of \cite{waki2013facial} for related statements.
\begin{lem} \label{lem:hp}
Let $\coneName \subseteq \mathcal{E}$ be a  convex cone and $\subspaceName \subseteq \mathcal{E}$ be an affine subspace for which $\subspaceName \cap \coneName$ is non-empty. The following
statements are equivalent.
\begin{enumerate}
\item $\subspaceName \cap \relint \coneName$ is empty.
\item There exists $s \in \coneName^* \setminus \coneName^{\perp}$ for which the hyperplane $s^{\perp}$ contains $\subspaceName$.
\end{enumerate}
  
\begin{proof}

To see (1) implies (2), note the main separation theorem (Theorem 11.3)  of Rockafellar \cite{rockafellar1997convex} states a hyperplane exists that \emph{properly} separates the sets $\subspaceName$ and $\coneName$ if the intersection of their
relative interiors is empty.  Using Theorem 11.7 of Rockafellar, we can additionally assume this hyperplane passes through the origin since $\coneName$  is a cone.  In
other words, if $\subspaceName \cap \relint \coneName$ is empty, there exists $s$ satisfying
\begin{align*} 
\begin{array}{cll}
\langle s, x \rangle \le 0  & \mbox{for all } & x \in  \subspaceName \\
\langle s, x \rangle \ge 0  & \mbox{for all }  & x \in  \coneName \\
\langle s, x \rangle \ne 0  & \mbox{for some } & x \in   \subspaceName \cup \coneName. 
\end{array}
\end{align*}
We will show that $\langle s, x \rangle = 0$ for all $x \in \subspaceName$, which will establish statement (2). Let $x_0$ denote a point in $\subspaceName \cap \coneName$ and let $\mathcal{T}$ be a subspace for which $\subspaceName = x_0 + \mathcal{T}$. 
Clearly, $\langle s, x_0\rangle = 0$. Since $\langle s, x_0\rangle$ vanishes, we must have that $\langle s, x \rangle  \le 0$  for all $x \in \mathcal{T}$. But $\mathcal{T}$ is a subspace,
therefore  $\langle s, -x \rangle \le 0$ also must hold. Thus,  $\langle s, x \rangle = 0$ for all $x \in \mathcal{T}$ and $s^{\perp}$ contains $\subspaceName$. 
Since $\langle s, x \rangle$ vanishes for all  $x\in \subspaceName$,  $\langle s, x \rangle \ne 0$ holds
for some $x \in \coneName$.  This establishes that $s$ is not in $\coneName^{\perp}$ and completes the proof.

To see that (2) implies (1), suppose $s \in \coneName^* \setminus \coneName^{\perp}$ exists and suppose for contradiction there exists an $x_0 \in \subspaceName \cap \relint \coneName$.  Since $\coneName^{\perp}$ is the  orthogonal complement of $\lin \coneName$, we can
 decompose $s$ as $s=s_1+s_2$, where $s_1 \in \lin \coneName$ and $s_2$ is in $\coneName^{\perp}$. Note that $s_1$ is also in $\coneName^*$,  $s_1$ is non-zero, and  $\langle s_1, x_0 \rangle = \langle s, x_0 \rangle = 0$. 
Since the affine hull of $\coneName$ equals the subspace $\lin \coneName$, $x_0$, being in the relative interior, is an interior point relative to $\lin \coneName$.  Hence, we must have that $x_0- \epsilon s_1$ is in $\coneName$ for some $\epsilon > 0$. This implies  
\[
\langle s_1, x_0- \epsilon s_1 \rangle = -  \epsilon ||s_1||^2 \ge 0,
\]
which cannot hold for any $\epsilon > 0$. Hence, no $x \in \subspaceName \cap \relint \coneName$  exists.

\end{proof}
\end{lem}
 The vector $s$  given by  statement (2) is called a \emph{reducing certificate} for $\subspaceName \cap \coneName$. 
Notice intersection with $s^{\perp}$  leaves  $\subspaceName \cap \coneName$ unchanged.  Letting $\faceName$ denote the face $\coneName \cap s^{\perp}$ therefore yields
the following equivalent  optimization problem:
\begin{align*} \label{opt:face}
	\begin{array}{cll} 
        \mbox{maximize } b^Ty   \mbox{ subject to }    c - Ay  \in \faceName.
	\end{array}
\end{align*}
Since faces of convex cones are also convex cones, this simplification can be repeated if one can find a reducing
certificate for $\subspaceName \cap \faceName$.  Indeed,  by taking $\faceName_{i+1} = \faceName_i \cap s^{\perp}_i$ for an $s_i
\in \faceName_i^* \setminus \faceName_i^{\perp}$  orthogonal to $\subspaceName$, one can find a chain of faces  $\faceName_i$
\[
\coneName = \faceName_0 \supset \faceName_{1} \supset \cdots \supset \faceName_{n-1} \supset \faceName_{n}
\]
that contain $\subspaceName \cap \coneName$. 
 An explicit algorithm for producing this chain is given in  \cite{pataki2013simple}, which we  reproduce in Algorithm~\ref{alg:fr}.

\begin{algorithm}[H] \caption{Facial reduction algorithm. Computes a sequences of faces $\faceName_i$ of the cone  $\coneName$  containing $\subspaceName \cap \coneName$,
where $\subspaceName$ is an affine subspace.} \label{alg:fr}
\Begin{
 Initialize: $\faceName_0 \leftarrow \coneName, i = 0$\\

    \Repeat{ $(\star)$ \rm is infeasible}{
\begin{enumerate}
\item \emph{Find reducing certificate,} i.e. solve the  feasibility problem    
\begin{align*}
\begin{array}{llc}
\mbox{{Find}} &  s_i \in \faceName_i^*   \setminus \faceName_i^{\perp} \\
 \mbox{{subject to  }} & s_i^{\perp} \mbox{ contains } \subspaceName 
\end{array} \;\;\;\;\;\;\; (\star)
\end{align*}
\item \emph{Compute new face,} i.e. set $\faceName_{i+1} = \faceName_i \cap s^{\perp}_i\;\;\; $ 
\item \emph{Increment counter $i$}
\end{enumerate}
      }
}
\end{algorithm}  
\noindent Note since the dimension of $\faceName_i$ drops at each iteration, Algorithm~\ref{alg:fr}  terminates after finitely many steps. Also note if the algorithm terminates after $n$ iterations, then $\subspaceName \cap \relint \faceName_n$ is non-empty by Lemma~\ref{lem:hp}. 

\paragraph{Finding reducing certificates}
To execute the facial reduction algorithm (Algorithm~\ref{alg:fr}), one must solve a feasibility problem $(\star)$ at each iteration to
 find $s_i \in \faceName^*_i \setminus \faceName_i^{\perp}$ orthogonal to $\subspaceName$.  It turns out this feasibility problem is also a conic optimization problem. To see this,  
recall the definition  of $\subspaceName$ from above (i.e. $\subspaceName :=\left\{ c - Ay : y \in \mathbb{R}^m \right\}$), let $A^* : \mathcal{E} \rightarrow \mathbb{R}^m$ denote the adjoint of $A$ and  pick $x_0$ in the relative interior of $\faceName_i$.  The solutions to $(\star)$ are (up to scaling) the solutions to:
\begin{align}
\begin{array}{lrc}
\mbox{{Find  }} &  s_i  \\
 \mbox{{subject to}}  & s_i \in \faceName_i^* ,  \langle s_i, x_0\rangle = 1 & \mbox{ (i.e. $s_i \in \faceName_i^* \setminus \faceName_i^{\perp}$)} \\
		      & A^* s_i = 0, \langle c,s_i \rangle = 0 & \mbox{ (i.e. $s_i^{\perp} \mbox{ contains } \subspaceName$ ).} 
\end{array} \label{eq:conePrg}
\end{align}
That $s^{\perp}_i$ contains $\subspaceName$  if and only if the second line of constraints holds can be shown using the standard identity  $(\range A)^{\perp} = \nullspace A^*$.  Correctness
of the first line of constraints  arises from the following corollary of Lemma~\ref{lem:hp}:
\begin{cor} \label{cor:relint}
Let $\coneName$ be a  convex cone, let $s$ be an element of $\coneName^*$, and
let $x$ be any element of $\relint \coneName$.  Then, $s \in \coneName^* \setminus \coneName^{\perp}$ if and only if $\langle s, x \rangle > 0$.
\begin{proof}
The if direction is obvious.  To see the other direction, suppose $s$ is in  $\coneName^* \setminus \coneName^{\perp}$ and $\langle s, x \rangle = 0$.   Applying Lemma~\ref{lem:hp}, this implies $\{ x \} \cap \relint \coneName$ is empty, a contradiction.
\end{proof}
\end{cor}

\paragraph{Discussion}
We make a few concluding remarks about the algorithm.  First, it terminates after finitely many steps, since the dimension of $\faceName_i$ drops at each iteration. 
Second, if the algorithm terminates after $n$ iterations, then $\subspaceName \cap \relint \faceName_n$ is non-empty, a simple consequence of Lemma~\ref{lem:hp}. In other words, a reformulation of the original problem over the face $\faceName_n$  is strictly feasible.

\begin{rem}
Throughout this section, we have assumed the given problem is feasible.  If the facial reduction algorithm (as presented) is applied to a problem that is infeasible,
it will identify  faces $\faceName_i$ for which the sets $\subspaceName \cap \faceName_i$ are also empty, leading to an equivalent problem that is also infeasible.  Though it is possible
to modify the algorithm to detect infeasibility (see, e.g., \cite{waki2013facial}),  we forgo this to simplify presentation.
\end{rem}

\subsection{Facial reduction of semidefinite programs} \label{subsec:sdpfr}

In this section, we develop a version of the facial reduction algorithm (Algorithm~\ref{alg:fr}) for semidefinite programs, i.e., we consider the case where the cone $\coneName = \mathbb{S}^{n}_{+}$
and the inner product space  $\mathcal{E} = \mathbb{S}^n$. 
This procedure, given   explicitly by Algorithm~\ref{alg:frDual}, represents each face $\faceName_i$ as a set of the form  $U_i \mathbb{S}^{d_i}_{+} U^T_i$ (with $d_i \le n$) for an appropriate
rectangular matrix $U_i$ (leveraging the description of faces given by Lemma~\ref{lem:facesofPSD}).   It finds reducing certificates $S_i \in \mathbb{S}^n$ by solving a semidefinite program over $\mathbb{S}_{+}^{d_i}$ and
 it computes a new face $\faceName_{i+1} := \faceName_i \cap S_i^{\perp}$ by finding a basis for the null space of particular matrix (related to the reducing certificate). 
 
 Algorithm~\ref{alg:frDual} applies to SDPs in the following form:
\begin{align*} 
\begin{array}{lll}
\mbox{\rm maximize} &  b^T y  \\
 \mbox{\rm subject to}  & C - \sum^{m}_{j=1}  y_j A_j   \in \mathbb{S}^n_{+},
\end{array} 
\end{align*}
where $C$ and $A_j$ are fixed symmetric matrices defining the following  affine subspace $\subspaceName$ of $\mathbb{S}^n$:
\[
\subspaceName := \left\{ C - \sum^{m}_{j=1}  y_j A_j  : y \in \mathbb{R}^m \right\}.
\]
 We now explain the basic steps of  Algorithm~\ref{alg:frDual} in more detail.

\begin{algorithm}
 Initialize: $U_0 = I_{n\times n}, d_0 = n$, $i = 0$ \\
    \Repeat{ $(\star)$ \rm is infeasible}{
\begin{enumerate}
 
\item \emph{Find reducing certificate $S_i$}, i.e. solve the  SDP
\begin{align*}
\begin{array}{lllccc}
\mbox{{Find}} &  S_i \in \mathbb{S}^{n}   \\
 \mbox{{subject to}} &  U_i^T   S_i U_i \in \mathbb{S}^d_{+},  U_i U_i^T \cdot S_i =1  & (i.e. \; S_i \in \faceName^{*}_i \setminus \faceName_i^{\perp}) \\
& C \cdot S_i  = 0,\;A_j \cdot S_i   = 0 \;\; \forall j \in \{1,\ldots,m\}   & (i.e. \; S_i^{\perp} \mbox{ contains } \subspaceName )
\end{array}(\star)
\end{align*}
\item \emph{Compute new face,} i.e. find  basis $B$ for $\nullspace U_i^T S_i U_i$,  set $U_{i+1}$ equal to \\ \qquad  $U_iB,$ and set $d_{i+1}$ equal to $\dim \nullspace  U_i^T S_i U_i$.
\item \emph{Increment counter $i$} 
\end{enumerate}
      }
\caption{Facial reduction algorithm for an SDP. Computes a sequences of faces $\faceName_i := U_i \mathbb{S}^{d_i}_{+} U^T_i$ of  $\mathbb{S}^n_{+}$  containing $\subspaceName \cap \mathbb{S}^n_{+}$,
where $\subspaceName :=  \left\{ C- \sum^m_{j=1} y_j A_j : y \in \mathbb{R}^m \right\}$. 
    }   \label{alg:frDual}
\end{algorithm} 

\paragraph{Step one: find reducing certificate}
At each iteration $i$,  Algorithm~\ref{alg:frDual} finds a reducing certificate $S_i \in \faceName^*_i \setminus \faceName_i^{\perp}$ for $\subspaceName \cap \faceName_i$, where $\faceName_i$ denotes the face $U_i \mathbb{S}^{d_i}_{+} U^T_i$.
This is done by solving   conic optimization  problem \eqref{eq:conePrg} specialized to the case $\coneName = \mathbb{S}^n_{+}$.
This specialization  appears as SDP $(\star)$, where we've used \eqref{eq:faceDualCnt} of Lemma~\ref{lem:facesofPSD} to describe $\faceName_i^*$ 
and the point $U_i U^T_i \in \relint \faceName_i$ to describe $\faceName_i^* \setminus \faceName_i^{\perp}$.

\paragraph{Step two: compute new face}
The second step of Algorithm~\ref{alg:frDual}  computes a new face by intersecting $\faceName_i$ with the subspace $S_i^{\perp}$. Computing this intersection 
can be done  using a matrix $B \in \mathbb{R}^{d \times r}$ with range equal to $\nullspace U_i^T S_i U_i$. Explicitly, we have that $\faceName_i \cap S_i^{\perp} = U_i B\mathbb{S}^r_{+} B^TU_i^T$, as  shown in the next lemma.
\begin{lem} \label{lem:intersect}
For $U \in \mathbb{R}^{n \times d}$, let $\faceName$ denote the set $U \mathbb{S}^d_{+} U^T$ and let $S \in \mathbb{S}^n$ and  $B \in \mathbb{R}^{d \times r}$ satisfy
\begin{align*}
 U^T S U \succeq 0,  \qquad  \range B = \nullspace U^T S U.
\end{align*}
 The following relationship holds:
\[
\faceName \cap S^{\perp} = U B \mathbb{S}^r_{+} B^T U^T.
\]

\begin{proof}
The containment $\supseteq$ is obvious.  To see the other containment, let $U X U^T$ be an  element of $\faceName \cap S^{\perp}$ for some $X \succeq 0$.  Taking inner product with $S$ yields
\[
 U X U^T \cdot S =  X  \cdot U^T S U  = 0.
\]
Since  $X \succeq 0$  and $U^T  S U \succeq 0$, the inner product $X  \cdot U^T S U$  vanishes if and only if $\range X$ is contained in $\nullspace U^T S U$  
(see, for example, Proposition 2.7.1 of \cite{pataki2000geometry}). In other words, $X$ is in the face $B \mathbb{S}^r_{+} B^T$ of $\mathbb{S}^d_{+}$, completing the proof.
\end{proof}
\end{lem}

\paragraph{Discussion}
We now make a few comments about Algorithm~\ref{alg:frDual}.  Variants of this algorithm arise by using different descriptions of $\faceName_i$ or by  using different
descriptions of the affine subspace $\subspaceName$.  If, for instance, one represents $\subspaceName$ as the set of $X$ solving the equations $A_j \cdot X = b_j$ for $j \in \{1,\ldots,m\}$, then
the set of $S_i$ orthogonal to $\subspaceName$ equals the set 
\begin{align} \label{eq:perpPrim}
 \left\{ \sum^m_{j=1} y_j A_j : y \in \mathbb{R}^m, b^Ty = 0 \right\}.
\end{align}
 Hence, to apply Algorithm~\ref{alg:frDual} to SDPs defined by equations $A_j \cdot X = b_j$, one simply replaces the constraints $C \cdot S_i = 0, A_j \cdot S_i = 0$ with membership in  \eqref{eq:perpPrim}.
We also note from Lemma~\ref{lem:facesofPSD} that one can  represent $\faceName_i$ and $\faceName_i^*$ using a sequence of invertible matrices $(U_i,V_i)$, 
which could be a more convenient description depending on implementation or the representation of $\subspaceName$.

\section{Our Approach} \label{sec:approach}
\subsection{Partial facial reduction} \label{sec:pfr}
Each iteration of the general facial reduction algorithm (Algorithm~\ref{alg:fr}) finds a reducing certificate by solving the feasibility problem $(\star)$.
Though the reducing certificate identifies a lower dimensional face, this benefit must be traded off with the cost of solving $(\star)$.
In this section, we propose a method for managing this trade-off.  Specifically, we describe a method for reducing the complexity of the feasibility problem $(\star)$ at the cost of only \emph{partially} simplifying the given conic optimization problem. 

Our method is as follows. Letting $\faceName_i$ denote the current face  at iteration $i$ of Algorithm~\ref{alg:fr}, we  approximate $\faceName_i$ with a user-specified convex cone $\faceName_{i,outer}$  that satisfies:
\begin{center}
\fbox{%
\parbox{0.5\linewidth}{ 
\begin{enumerate}
 \item $\faceName_i \subseteq \faceName_{i,outer}$ \mbox{ (which implies  $\faceName_{i,outer}^* \subseteq \faceName_i^*$})
 \item $\lin \faceName_i = \lin \faceName_{i,outer}$ (i.e. $\faceName^{\perp}_{i,outer} = \faceName_i^{\perp}$)
 \item $\faceName_{i,outer}^*$ has low search complexity,
\end{enumerate}
}
}
\end{center}
where the first two conditions ensure that $\faceName_{i,outer}^* \setminus  \faceName^{\perp}_{i,outer}$ is a subset of $\faceName_{i}^* \setminus \faceName_{i}^{\perp}$.
Using the approximation $\faceName_{i,outer}$, we then modify the feasibility problem $(\star)$ to search over this subset:
\begin{center}
\fbox{%
\parbox{0.6\linewidth}{ 
\begin{align*}
\begin{array}{llc}
\mbox{{Find}} &  s_i \in \faceName_{i,outer}^* \setminus  \faceName^{\perp}_{i,outer} \subseteq \faceName_{i}^* \setminus \faceName_{i}^{\perp} \\
 \mbox{{subject to  }}   & s_i^{\perp} \mbox{ contains } \subspaceName.   
\end{array} \;\;\;\;\;\;\; (\star)
\end{align*}
}
}
\end{center}
By construction, a solution $s_i \in \faceName_{i,outer}^* \setminus  \faceName^{\perp}_{i,outer}$ to the modified feasibility problem is a solution to the original; hence, $\faceName_i \cap s^{\perp}_i$ 
is a face  of $\faceName_i$ (and the cone $\coneName$) containing    $\subspaceName \cap \coneName$. Further, the approximation $\faceName_{i,outer}$ can be chosen such that the search complexity of $\faceName_{i,outer}^*$ matches desired pre-processing effort. In other words, the algorithm correctly identifies a face   at cost specified by the user. 

\subsubsection{Existence of reducing certificates} \label{sec:redexist}

\begin{figure}
\centering
\includegraphics[scale=.6]{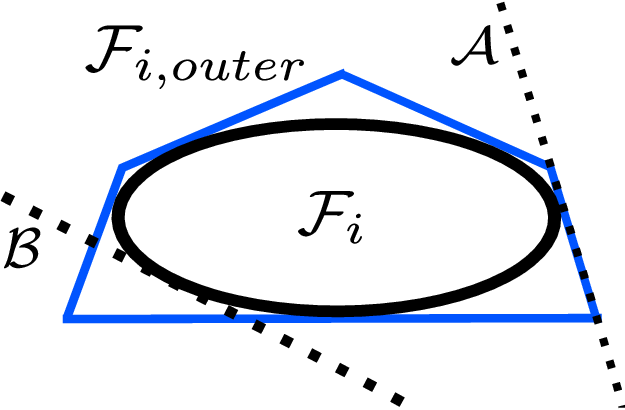}
\caption{Illustrates when the facial reduction algorithm successfully finds a face when modified to use the approximation $\faceName_{i,outer}$.   If the feasible set is  contained in $\subspaceName \cap  \faceName_i$,
the procedure succeeds:  $\subspaceName \cap \relint \faceName_{i,outer}$ is empty. In contrast, if the feasible set is  contained in $\mathcal{B} \cap  \faceName_i$,  
the procedure fails:  $\mathcal{B} \cap \relint \faceName_{i,outer}$ is non-empty.  
 }
\label{fig:coneRelax}
\end{figure}

Because we have introduced  the approximation $\faceName_{i,outer}$, the algorithm may not find a reducing certificate (and hence a lower dimensional face) even if $\subspaceName \cap \relint \faceName_i$ is empty; hence,
the algorithm may not find  a face of minimal dimension. This leads to the following question: when  will the modified feasibility problem $(\star)$ have a solution? Since we have chosen $\faceName_{i,outer}$ to be a convex
cone, we can use Lemma~\ref{lem:hp} to answer this question.  Under
the assumption that $\subspaceName \cap \faceName_{i,outer}$ is
non-empty, this lemma states feasibility of $(\star)$ is now
\emph{equivalent} to emptiness of $\subspaceName \cap \relint
\faceName_{i,outer}$. In other words, the modified feasibility problem $(\star)$ has a solution
 if and only if a \emph{relaxation} of the problem of interest is not strictly feasible.
Figure~\ref{fig:coneRelax} illustrates a situation when this condition holds and when it fails
for two different subspaces.

\subsubsection{Approximating faces of $\mathbb{S}^n_{+}$}
To apply this idea to SDP, and to  modify the SDP facial reduction algorithm (Algorithm~\ref{alg:frDual}), we need a way of approximating faces of $\mathbb{S}^n_{+}$. To see how this can be done, let $\faceName$ denote a face $U \mathbb{S}^d_{+} U^T$ of $\mathbb{S}^n_{+}$
for some $U \in \mathbb{R}^{n \times d}$.  An approximation $\faceName_{outer}$  can be defined using an approximation $\hat{\mathbb{S}}^{d}_{+}$  of $\mathbb{S}^d_{+}$.
Moreover, the search complexity of $\faceName^*_{outer}$ depends on the search complexity of $\hat{\mathbb{S}}^{d}_{+}$.
Consider the following (whose proof is straight-forward and omitted):
\begin{lem} \label{lem:psdFaceAppx}
Let  $\hat{\mathbb{S}}^{d}_{+} \subseteq \mathbb{S}^d$ be a convex cone containing $\mathbb{S}^{d}_{+}$. For   $U \in \mathbb{R}^{n \times d}$, let $\faceName_{outer}$ and $\faceName$ denote the sets $U \hat{\mathbb{S}}^{d}_{+} U^T$
and $U {\mathbb{S}}^{d}_{+} U^T$, respectively. The following statements are true.
\begin{enumerate}
 \item  $\faceName \subseteq \faceName_{outer}$.
\item $\lin \faceName = \lin \faceName_{outer}$
\item $\faceName^*_{outer} = \left\{ X \in \mathbb{S}^n : U ^T X U \in (\hat{\mathbb{S}}^{d}_{+})^* \right\}$.
\end{enumerate}

\end{lem}
\noindent Based on this lemma, we conclude to modify  Algorithm~\ref{alg:frDual},
 it suffices to  replace the PSD constraint of  SDP $(\star)$   with membership in  $(\hat{\mathbb{S}}^{d}_{+})^*$,
where $\hat{\mathbb{S}}^{d}_{+}$  is a cone outer-approximating  $\mathbb{S}^d_{+}$. Example
approximations are explored in the next section.
 
\subsection{Approximations of $\mathbb{S}^d_{+}$} \label{sec:approx}
In this section, we explore an  outer approximation $\outerName$ of  $\mathbb{S}^d_{+}$  parametrized by a set $\setName$
of $d \times k$ rectangular matrices. The parametrization is chosen such that the dual cone $\approxName$ equals the \emph{Minkowski sum} of faces $W_i \mathbb{S}^{k}_{+} W_i^T$ of $\mathbb{S}^d_{+}$ 
for $W_i \in \setName$.  It is defined below:
\begin{lem}
For a set $\setName := \left\{ W_1, W_2,\ldots, W_{|\setName|} \right\}$ of $d \times k$ matrices, let $\outerName$ denote the following convex cone:
\begin{align*}
\outerName := \left\{  X \in  \mathbb{S}^d : W_i^T X W_i \in \mathbb{S}^k_{+} \;\;\;i = 1, \ldots, |\setName| \right\}.
\end{align*}
The dual cone $\approxName$ satisfies
\begin{align} \label{def:coneProdInner}
 \approxName = \left\{ \sum^{|\setName|}_{i=1}  W_i X_i W_i^T : X_i \in \mathbb{S}_{+}^k  \right\},
\end{align}
and the following inclusions hold:
\[
\approxName \subseteq \mathbb{S}^d_{+} \subseteq \outerName.
\]
\begin{proof}
The inclusions are obvious from the definitions of $\approxName$ and $\outerName$ (as is the fact that  $\outerName$ is a convex cone).  It remains to show correctness of  \eqref{def:coneProdInner}. To show this,
let $\mathcal{T}$ denote the set on the right-hand side of \eqref{def:coneProdInner}.  It is easy to check that  $\mathcal{T}^*= \outerName$, which implies $\mathcal{T}^{**} = \approxName$.
Since $\mathcal{T}$ is a convex cone (as is easily checked), $\mathcal{T}^{**}$ equals the closure of $\mathcal{T}$. The result therefore follows by showing $\mathcal{T}$ is closed. To see this, note that $\mathcal{T}$ equals the Minkowski sum of closed cones $W_i \mathbb{S}^k_{+} W_i^T$.
For  matrices $Z_i \in W_i \mathbb{S}^k_{+} W_i^T$, we have that $\sum^{|\setName|}_{i=1} Z_i=0$
only if $Z_i= 0$ for each $i$.  This shows that $\sum^{|\setName|}_{i =1}Z_i = 0$ only if $Z_i$ is in the lineality space of $W_i \mathbb{S}^k_{+} W_i^T$.
 Direct application of the closedness criteria Corollary 9.1.3 of Rockafellar \cite{rockafellar1997convex} shows $\mathcal{T}$ is closed.
\end{proof}

\end{lem}

Since the modification to the SDP facial reduction algorithm (Algorithm~\ref{alg:frDual}) will involve searching over $\approxName$ (as indicated by Lemma~\ref{lem:psdFaceAppx}), we will investigate $\outerName$ by studying the dual cone $\approxName$.
We first make a few comments regarding the search complexity of $\approxName$ for different choices of $\setName$. 
 Note when $k=1$, each $W_j$ in $\setName$ is a vector and 
$\approxName$  is the conic hull of a finite set of rank one matrices.  In other words, $\approxName$ is \emph{polyhedral} and  can be described  by \emph{linear programming}.  When $k=2$, 
the set $\approxName$ is defined by $2 \times 2$ semidefinite constraints and can hence be described by \emph{second-order cone programming} (SOCP). This follows since each $X_i \in \mathbb{S}_{+}^2$ can be expressed using scalars $a,b,c$ constrained as follows:
\begin{align} \label{eq:socp} 
X_i = \left(  \begin{array}{cc}
    a + b & c \\
	c &  a-b
\end{array} \right) \succeq 0 \;\;\;\;\; \Leftrightarrow  \;\;\;\;\; a \ge 0 \;\;\; \mbox{ and }  a^2 \ge b^2 + c^2.
\end{align} 
Example choices for $\approxName$ are now given.  As we will see, well-studied approximations of $\mathbb{S}^d_{+}$  can be expressed as sets of the form  $\approxName$.  

\subsubsection{Examples}
Example choices for $\approxName$ are given in Table~\ref{tab:manyCones} along with the cardinality of the set $\setName$ that yields each entry.
\begin{table}
\begin{center}
    \begin{tabular}{ | c| c | c   | c |}
    \hline
   $\outerName$ &  $\approxName$ & Search &  $|\setName|$  \\ \hline 
   $X_{ii} \ge 0$ &   Non-negative diagonal $(\D^d)$ & LP  & $\bigO(d)$ \\ \hline 
    $X_{ii} \ge 0$, $X_{jj}+X_{ii} \pm 2 X_{ij}\ge 0$ & Diagonally-dominant $(\DD^d)$& LP    & $\bigO(d^2)$  \\  \hline 
    $2\times2$ principal  sub-matrices psd & Scaled diagonally-dominant $(\SDD^d)$& SOCP &  $\bigO(d^2)$ \\ \hline 
    $k\times k$ principal  sub-matrices psd  & Factor width-$k$ ($\FW^d_{k}$)  &   SDP &  $\bigO({d  \choose  k} )$ \\ \hline 
    \end{tabular}
\end{center}
   \caption{Example outer and inner approximations of $\mathbb{S}^d_{+}$, the search algorithm for $\approxName$,  and the cardinality of the set $\setName$.} \label{tab:manyCones}
\end{table}
 Included are $d \times d$ non-negative diagonal matrices $\D^d$,  diagonally-dominant matrices $\DD^d$, scaled diagonally-dominant matrices  
$\SDD^d$ as well as matrices $\FW^d_k$  with \emph{factor-width}  \cite{boman2005factor} bounded by $k$. These sets satisfy
\[
 \D^d = \FW^d_1 \subseteq \DD^d \subseteq \SDD^d = \FW^d_2  \subseteq \FW^d_{3} \subseteq \cdots  \subseteq \FW^d_{d} = \mathbb{S}^{d}_{+},
\]
and the sets $\D^d$ and $\DD^d$ are polyhedral. Details on each entry follow.

\subsubsection{Non-negative diagonal matrices ($\D^d$) }
A simple choice for $\approxName \subseteq \mathbb{S}^d_{+}$ is the set of non-negative diagonal matrices:
\begin{align*}
\D^d := \left\{ X \in \mathbb{S}^d :  X_{ii} \ge 0,  \;\; X_{ij} = 0 \;\; \forall i \ne j \right\}.
\end{align*}
The set $\D^d$  contains non-negative combinations of matrices $w_iw_i^T$, where $w_i$ is a permutation of $(1,0,\ldots,0,0)^T$.
In other words, the set $\D^d$ corresponds to the set $\approxName$ if we take 
\[
 \setName = \left\{ (1,0,\ldots,0,0)^T, (0,1,\ldots,0,0)^T,\ldots, (0,0,\ldots,0,1)^T \right\}.
\]

\subsubsection{Diagonally-dominant matrices ($\DD^d$) } \label{subsec:diagDom}
Another well studied choice for $\approxName$  is cone of  symmetric diagonally-dominant matrices 
with non-negative diagonal entries \cite{barker1975cones}:
\[
\DD^d := \left\{ X \in \mathbb{S}^d :  X_{ii} \ge \sum_{j \ne i} |X_{ij}| \right\}.
\]
This set is  polyhedral.  The \emph{extreme rays} of $\DD^d$   are matrices of the form $w_i w_i^T$, where $w_i$ is any permutation of
\[
(1,0,0,\ldots,0)^T,  (1,1,0,\ldots,0)^T,  \mbox{ or }  (1,-1,0,\ldots,0)^T.
\]
 Taking $\setName$ equal to the set of all such permutations gives $\approxName = \DD^d$. This representation makes the inclusion $\DD^d \subseteq \mathbb{S}^d_{+}$ obvious. We also see that $\DD^d$ contains $\D^d$.

\subsubsection{Scaled diagonally-dominant  matrices ($\SDD^d$)}
A non-polyhedral generalization of $\DD^d$ is the set of \emph{scaled diagonally-dominant matrices} $\SDD^d$.   This set equals all matrices obtained by pre- and post-multiplying  diagonally-dominant matrices
by   diagonal matrices with strictly positive diagonal entries:
\[
\SDD^d := \left\{ D T D : D \in\D^d, D_{ii} > 0, T \in \DD^d \right\}.
\]
The set $\SDD^d$  can be equivalently defined as the set of matrices that equal the sum of PSD matrices non-zero only on a $2\times 2$ principal sub-matrix (Theorem~9 of~\cite{boman2005factor}). 
 As an explicit example, we have that $\SDD^3$ are all matrices $X$ of the form
\begin{align*}
X = \underbrace {  \left( \begin{array}{ccc} a_{11} & a_{12} & 0 \\  a_{12} & a_{22}  & 0 \\ 0 & 0  & 0 \end{array} \right)}_{X_1} + 
      \underbrace{ \left( \begin{array}{ccc} b_{11} & 0 & b_{13} \\  0 & 0  & 0 \\ b_{13} & 0  & b_{33} \end{array} \right)}_{X_2} + 
        \underbrace{  \left( \begin{array}{ccc} 0 & 0 & 0 \\  0 & c_{22}  & c_{23} \\ 0 & c_{23}  & c_{33} \end{array} \right)}_{X_3},
\end{align*}
where $a_{ij}$, $b_{ij}$, and $c_{ij}$ are scalars chosen such that $X_1, X_2$ and $X_3$ are PSD. In general,  $\SDD^d$  equals $\approxName$  when $\setName$ 
 equals the set of  $d \times 2$ matrices $W$ for which $W^TXW$ returns a $2 \times 2$ principal sub-matrix of $X$. For $\SDD^3$, we have
\[
\SDD^3 = \mathcal{C}( \{ W_1, W_2,W_3 \} )^* =  \left\{  \sum^{3} _{i=1} W_i X_i W_i^T : X_i \in \mathbb{S}^2_+ \right\},
\]
where
\begin{align*}
\begin{array}{ccc}
W_1 = \left( \begin{array}{cc} 1 & 0 \\ 0 & 1 \\ 0 & 0 \end{array} \right) &  W_2 = \left( \begin{array}{cc} 1 & 0 \\ 0 & 0 \\ 0 & 1 \end{array} \right) & W_3 = \left( \begin{array}{cc} 0  & 0 \\ 1 & 0 \\ 0 &1 \end{array} \right).
\end{array}
\end{align*}
Also note from \eqref{eq:socp}  that $\SDD^d$ can be represented using second-order cone constraints.  This latter fact is used in recent work of Ahmadi and Majumdar \cite{AhmadiM14} to define an SOCP-based 
method for testing polynomial non-negativity. (A similar LP-based method is also presented in  \cite{AhmadiM14} that incorporates $\DD^d$.)

\paragraph{The kernels of $\SDD$ matrices}
The kernel of a scaled diagonally-dominant matrix has a structured  basis of vectors with \emph{disjoint support}, where the support of a vector $u \in \mathbb{R}^n$ is the set of indices $i$ for which $u_i \ne 0$.   This follows because, up-to permutation, a scaled diagonally-dominant is block-diagonal, where each block is either positive definite, equals the zero matrix, or has co-rank one (i.e., has a one dimensional kernel), as shown in \cite{chen2004combinatorial}. In Section~\ref{sec:simpleEx}, we use this result to show that reduced SDPs can be formulated without damaging sparsity if $\SDD$-approximations  are used (which, since $\D^d \subseteq \DD^d \subseteq \SDD^d$, shows sparsity is not damaged when diagonally-dominant or diagonal approximations are used). The following proposition summarizes relevant results of \citet{chen2004combinatorial}.  We include an elementary---and different---proof for completeness.
\begin{prop} \label{prop:kernelsdd}
Suppose $X \in \mathbb{S}_+^d$ is scaled-diagonally dominant.  Then, there is a permutation matrix $P \in \mathbb{R}^{d \times d}$ for which
\begin{align} \label{proof:graph}
P X P^T = \left( \begin{array}{cccc} 
X_1 & 0 & \cdots & 0 \\
0 & X_2 & \cdots & \vdots \\
\vdots & \vdots & \ddots & 0 \\
0  & 0 & 0 & X_M \\
\end{array} \right),
\end{align}
where,  for all $m \in \{1,\ldots,M\}$, the matrix $X_m \in \mathbb{S}_+^{d_m}$ is either positive definite, a matrix of all zeros, or has co-rank one. Moreover, when $X$ has co-rank $r$, there is a matrix $U \in \mathbb{R}^{d \times r}$ whose columns have disjoint support and span the nullspace of $X$.

\begin{proof}
For $X \in \mathbb{S}^d$, let $G_X := ([d], E)$ denote the  graph with node set $[d] := \{1,\ldots,d\}$, where $\{i,j\}$ is in the edge set $E$ if and only if $X_{ij} \ne 0$.   Clearly there is a permutation matrix $P$ that block-diagonalizes $X$ as in \eqref{proof:graph}, defined in the obvious way by the connected-components of $G_X$.  

Now suppose $P$ in \eqref{proof:graph} equals this permutation. That $X_m$ has the claimed properties is immediate when $d_m \le 2$.  Now suppose $d_m > 2$ and that $X_m$ is non-zero and not positive definite. Also, define the graph  $G_{X_m} = ([d_m],E_m)$, where $\{i,j\}$ is in the edge set $E_m$ if and only if $[X_m]_{ij} \ne 0$ and observe $G_{X_m}$ is connected (and, indeed, isomorphic to a connected component of $G_X$ defined above.)

 We first claim all components of $v  \in \nullspace X_m \setminus \{0\}$ are non-zero. To begin, pick $i \in [d_m]$ such that $v_i$ is non-zero.  For arbitrary $t \in [d_m] \setminus i$, there is a  path $T \subseteq E_m$ from $i$ to $t$ for which
 \[
 X_m = \bar{X} + \sum_{ \{r,s\} \in T } (e_r,e_s) X_{rs} (e_r,e_s)^T ,
 \]
where all entries of $X_{rs} \in \mathbb{S}^2_{+}$ are non-zero and $\bar{X}$ is positive semidefinite. Picking the first edge $\{i,j\} \in T$, we conclude that $X_{ij} (e_i,e_j)^Tv = X_{ij}(v_i,v_j)^T = 0$.
For the sake of contradiction, suppose $v_j = e_j^Tv = 0$. Then, $(v_i,0)^T$ is in the kernel of $X_{ij}$, showing a diagonal entry of $X_{ij}$ is zero (since $v_i \ne 0$), contradicting the fact all entries of $X_{ij}$ are non-zero. Hence, $v_j \ne 0$.  Repeating this argument  using the next edge $\{j,k\}$ in $T$ shows $v_k \ne 0$.  Repeating for all edges in $T$  shows $v_t \ne 0$. Since $t$ was arbitrary, all components of $v$ are non-zero.

Now pick another non-zero $w \in \nullspace X_m$ and consider the consecutive edges $\{i,j\}$ and $\{j,k\}$ in the path $T$. Then, for scalars $\lambda$ and $\gamma$,
\[
(e_i,e_j)^T w = \lambda (e_i,e_j)^Tv, \;\; (e_j,e_k)^Tw = \gamma (e_j,e_k)^Tv,
\]
 otherwise the non-zero matrices $X_{ij} \in \mathbb{S}^2_{+}$ and $X_{jk} \in \mathbb{S}^2_{+}$ have two-dimensional kernels, and are therefore the zero matrix, a contradiction. But since  $v_j$ and $w_j$ are non-zero, we also have $\lambda = \gamma$.  Since any $s,t \in [d_m]$ are connected by a path,  we conclude $w = \lambda v$.

Existence of $U$ is immediate, given that the  kernel of $X_m$ has a basis of the form $\{e_1,\ldots,e_{d_m}\}$, $\{0\}$, or $\{v\}$.

\end{proof}

\end{prop}

\subsubsection{Factor-width-k matrices}
A generalization of $\SDD^d$ (and diagonal matrices $\D^d$) arises from notion of \emph{factor-width} \cite{boman2005factor}.
The factor-width of a matrix $X$ is the smallest integer $k$ for which  $X$ can be written as the sum of PSD matrices that are non-zero only on a single $k\times k$ principal sub-matrix.  

Letting $\FW^d_k$   denote the set of $d \times d$  matrices of factor-width no greater than $k$, we have that $\SDD^d =  \FW^d_2$ and
$\D^d =  \FW^d_1$.  To represent $\FW^d_k$ as a cone of the form $\approxName$, we set $\setName$ to be the set of  $d \times k$  matrices $W_j$ for which $W_j^TXW_j$ returns a $k \times k$ principal sub-matrix of $X$.
Note that there are ${d \choose k}$ such matrices, so a complete parametrization of  $\FW^d_k$ is not always practical using this representation.  Also note  $\FW^d_k$ equals $\mathbb{S}^{d}_{+}$ when $k=d$.

\subsubsection{Corresponding outer approximations}

We briefly discuss the outer approximation $\outerName$ corresponding to the discussed examples for $\approxName$. To summarize, if $\approxName$ is the cone of non-negative diagonal matrices $\D^d$, then
$\outerName$ is the cone of matrices whose diagonal entries are non-negative. If $\approxName$ is the cone of non-negative diagonally-dominant matrices $\DD^d$, then $\outerName$ is the set of matrices $X$
for which $w_i^T X w_i \ge0$, where $w_iw_i^T$ is an extreme ray of $\DD^d$ (given  in Section~\ref{subsec:diagDom}).  If $\approxName$ is the set of scaled diagonally-dominant matrices $\SDD^d$, 
then $\outerName$  is the set of matrices with positive semidefinite $2\times2$ principal sub-matrices. 
Finally, if $\approxName$ equals $\FW^d_k$, the set of matrices with factor-width bounded by $k$, then $\outerName$  is the set of matrices with positive semidefinite $k\times k$ principal sub-matrices.  We see as $\approxName$ grows larger, the constraints 
defining $\outerName$ become more restrictive, equaling a positive semidefinite constraint when $d=k$.

\subsection{Finding faces of minimal dimension/rank maximizing reducing certificates}
Suppose $\faceName_i := U_i \mathbb{S}^{d_i}_{+}  U_i^T$ is the current face at iteration $i$ of  the SDP facial reduction algorithm (Algorithm~\ref{alg:frDual}). Further suppose
 $\faceName_{i,outer}  := U_i \mathcal{C}(\setName_i) U^T_i$ approximates $\faceName_i$ per the discussion in Section~\ref{sec:pfr} (for some specified set of rectangular matrices $\setName_i$). The following question is natural: how can one find a reducing
certificate $S_i$ that   minimizes the dimension of  the  face $\faceName_{i+1} := \faceName_i \cap S_i^{\perp}$ when $S_i$ is constrained to  $\faceName^*_{i,outer} \setminus \faceName_{i,outer}^{\perp}$?
Using Lemma~\ref{lem:intersect}, it is easy to see this problem is solved by  finding a solution to 
\begin{align*}
\begin{array}{lllccc}
\mbox{{Find}} &  S_i \in \mathbb{S}^{n}   \\
 \mbox{{subject to}} &  U_i^T   S_i U_i \in \mathcal{C}(\setName_i)^*   & (i.e. \; S_i \in \faceName_{i,outer}^{*} ) \\
& C \cdot S_i  = 0,\;A_j \cdot S_i   = 0 \;\; \forall j \in \{1,\ldots,m\}   & (i.e. \; S_i^{\perp} \mbox{ contains } \subspaceName )
\end{array} 
\end{align*}
that  maximizes the rank of $U_i^T S_i U_i$.  In this section, we give a method for finding solutions of this type. 

To ease notation, we drop the subscript $i$ and also consider a more general question: how does one find maximum rank matrices in the set $\perpName \cap \approxName$, when $\perpName$ is an arbitrary linear subspace? 
(In the above, $\perpName$ is the subspace $\left\{  U_i^T S_i U_i : C \cdot S_i = 0, A_j \cdot S_i = 0 \right\}$.)  An answer to this question  arises  from the next two lemmas.

\begin{lem} \label{lem:rankX_i} 
Let $\perpName$ be a subspace of $\mathbb{S}^d$. If $X^\star:=\sum^{|\setName|}_{i=1} W_i X^\star_i W^T_i$  maximizes   $\sum^{|\setName|}_{i=1} \rank X_i$ over $\perpName \cap \approxName$, then    $X^\star$  maximizes $\rank X$ over $\perpName \cap \approxName$.
\begin{proof}

We will  argue the kernel of $X^\star$ is contained in the kernel of any  $X \in \perpName \cap \approxName$, which immediately implies $\rank X^\star \ge \rank X$. 

To begin, we first argue for any $X = \sum^{|\setName|}_{i=1} W_i X_i W^T_i \in  \perpName \cap \approxName$ that $\nullspace X^\star_i \subseteq \nullspace X_i$ for all $i \in \{1,\ldots,|\setName|\}$. To see this, first  note  that for any $X  \in  \perpName \cap \approxName$
the matrix
\[
 X^\star + X =  \sum^{|\setName|}_{i=1} W_i ( X^\star_i + X_i )W^T_i 
\]
is also in  $\perpName \cap \approxName$ and  satisfies $\rank (X^\star_i + X_i) \ge \rank X^\star_i$.  Now suppose for some $d \in \{1,\ldots,|\setName|\}$ that  $\nullspace X^\star_d \not\subseteq \nullspace X_d$. This implies that  $\nullspace (X^\star_d + X_d) =  \nullspace X^\star_d \cap \nullspace X_d  \subset \nullspace X^\star_d$ which in turn implies $\rank (X^\star_d + X_d ) > \rank X^\star_d$.
But this contradicts our assumption that  $X^\star$ maximizes $\sum_i \rank X_i$.  Hence, $\nullspace X^\star_i \subseteq \nullspace X_i$ for all $i \in \{1,\ldots,|\setName|\}$.

Now suppose an  $X \in \perpName \cap \approxName$ exists for which  $X^\star w = 0$ but $Xw \ne 0$ for some $w$.  Since $X w = 0$ if and only if $X_i W_i^T w = 0$ for all $i$, we must have for some $d$ that $W_d^Tw$ is in the kernel of $X^\star_d$ but not in the kernel of $X_d$. But we have already established that $\nullspace X^\star_d \subseteq \nullspace X_d$.  Hence, $w$ cannot exist. We therefore have that $\nullspace X^\star \subseteq \nullspace X$ for any $X\in \perpName \cap \approxName$, which completes the proof. 
\end{proof}
\end{lem}

We can use this  condition to formulate an SDP whose optimal solutions yield maximum rank matrices of $\perpName \cap \approxName$.
To maximize $\sum^{|\setName|}_{i=1} \rank X_i$, we introduce matrices $T_i$ constrained such that their traces $\trace T_i$ lower bound
$\rank X_i$.  We then optimize the sum of their traces.

\begin{lem}  \label{lem:sumOfConesMaxRank}
Let $\perpName$ be a subspace of $\mathbb{S}^d$.  A matrix $X$ maximizing $\sum^{|\setName|}_{i=1} \rank X_i$  over $\perpName \cap \approxName$  is given by any optimal solution $(X,X_i,T_i)$  to the following SDP:
\begin{align} \label{sdp:sumOfConesMaxRank}
\begin{array}{lll}
 \mbox{\rm maximize}  &  \sum^{|\setName|}_{i=1} \trace{T_i} \  \\
 \mbox{\rm subject to } & X \in \perpName,\\
&X = \sum^{|\setName|}_{i=1} W_i X_i W^T_i   & \mbox{ \rm  i.e. } X \in \approxName \\
& X_i \succeq T_i  & \forall i\in\{1, \ldots, |\setName|\}\\
&  I \succeq  T_i \succeq 0 &  \forall i\in\{ 1, \ldots, |\setName|\}.
\end{array} 
\end{align}

\begin{proof}
Let $r_{\max}$ equal the maximum of  $\sum^{|\setName|}_{i=1} \rank X_i$ over the  set of feasible $X_i$.  We will
show at optimality $\sum^{|\setName|}_{i=1} \rank X_i = r_{\max}$. 

To begin, the constraint $I \succeq T_i \succeq 0$ implies the eigenvalues of $T_i$ are  less than one. Hence,  $\rank T_i  \ge \trace T_i $.  
Since  $X_i \succeq T_i$ , we also have $\rank{X_i} \ge \rank{T_i}$.  Thus,  any feasible $(X_i,T_i)$ pair satisfies
\begin{align} \label{ineq:maxrank}
r_{\max} \ge \sum^{|\setName|}_{i=1} \rank{X_i}  \ge   \sum^{|\setName|}_{i=1} \rank{T_i}  \ge  \sum^{|\setName|}_{i=1} \trace T_i.
\end{align}

Now note for any feasible $(X,X_i)$ we can pick $\alpha > 0$  and construct a feasible point $(\alpha X, \alpha X_i, \hat{T}_i)$  that satisfies $\sum^{|\setName|}_{i=1} \trace{\hat{T_i}} =  \sum^{|\setName|}_{i=1}\rank{X_i}$; if $X_i$ has  eigen-decomposition $\sum_j \lambda_j u_j u_j^T$ for $\lambda_j > 0$, simply take  $\hat{T_i} = \sum_j u_j u_j^T$ and $\alpha$ equal to
\[
\max \bigcup_i \left\{ \frac{1}{\lambda} : \lambda \mbox{ is a  positive eigenvalue of $X_i$ } \right\}.
\]  
Hence, some feasible point $(\hat{X},\hat{X_i},\hat{T_i})$ satisfies  $\sum^{|\setName|}_{i=1} \trace{\hat{T_i}} = r_{\max}$. Therefore, the optimal $(X,X_i,T_i)$ satisfies 
\[
\sum^{|\setName|}_{i=1} \trace T_i  \ge  r_{\max}.
\]
Combining this inequality with \eqref{ineq:maxrank} yields that at optimality
\[
\sum^{|\setName|}_{i=1} \trace T_i   = \sum^{|\setName|}_{i=1}\rank{X_i}  = r_{\max},
\]
which completes the proof.
\end{proof}
\end{lem}
\noindent Combining the previous two lemmas shows how to maximize rank over $\perpName \cap \approxName$:
\begin{cor} \label{cor:maxRank}
 A matrix $X \in \perpName \cap \approxName$ of maximum rank is given by any optimal solution
 $(X,X_i,T_i)$   to the SDP \eqref{sdp:sumOfConesMaxRank}. 
\end{cor}

\paragraph{Maximum rank solutions for polyhedral approximations}
 We next illustrate how the search for maximum rank solutions simplifies when $\approxName$ is polyhedral.
Recall if $\setName$ is a set of vectors, i.e. if $k=1$, then $\approxName$ is the conic hull of a finite set of rank one matrices. 
In other words,  $\approxName$  is the set of matrices of the form $\sum^{|\setName|}_{i=1} \lambda_i w_i w_i^T$ 
for $\lambda_i \ge 0$ and $w_i \in \setName$.  In this case,   SDP~\eqref{sdp:sumOfConesMaxRank}  simplifies into the following linear program.
\begin{cor} \label{lem:lpPolyMaxRank}
 A matrix $X \in \perpName \cap \approxName$ of maximum rank is given by any  optimal  solution $(X, \lambda, t)$ to the following LP:
\begin{align}
\begin{array}{lll}
\emph{ \mbox{maximize }} &   \sum^{|\setName|}_{i=1} t_i\  \\
\emph{ \mbox{subject to  }} & X \in \perpName \\
&X = \sum^{|\setName|}_{i=1} \lambda_i w_i w^T_i ,  & \mbox{ \rm i.e. } X \in \approxName \\
& \lambda_i \ge t_i   & \forall i\in\{ 1, \ldots ,|\setName| \}\\
&1 \ge t_i \ge 0 &  \forall i\in\{ 1, \ldots, |\setName| \}.
\end{array} \label{lp:polyMaxRank}
\end{align}
\end{cor}

\paragraph{An alternative approach}
We  mention an alternative to SDP~\eqref{sdp:sumOfConesMaxRank} for maximizing $\sum^{|\setName|}_{i=1} \rank X_i$.  Notice that membership in $\approxName$
can be expressed using a semidefinite constraint on a block-diagonal matrix,  where maximizing the rank of this matrix is equivalent to maximizing $\sum^{|\setName|}_{i=1} \rank X_i$.
We conclude  $\sum^{|\setName|}_{i=1} \rank X_i$ is maximized  by finding a maximum  rank  solution to  a particular block-diagonal SDP, which can be done using interior point methods (since solutions
in the relative interior of the feasible set are solutions of maximum rank). Note, however, that
this alternative approach does not permit use of the simplex method when $\approxName$ is polyhedral, since the simplex method produces solutions on the boundary of the feasible set.  In contrast, the simplex
method can be used to solve LP \eqref{lp:polyMaxRank}, the specialization of SDP~\eqref{sdp:sumOfConesMaxRank}  to   polyhedral $\approxName$.

\subsection{Strictly feasible formulations}

It is well-known   the problem $(\star)$ of Algorithm~\ref{alg:frDual} is not
strictly feasible.  While this does not preclude solution of $(\star)$ via
the self-dual embedding technique (e.g., \cite{de1997initialization},\cite{ye1994nl}),
it does preclude solution via standard barrier methods  and may cause numerical difficulties. It is therefore natural to ask if our   formulation  for maximum rank reducing certificates  is always strictly feasible.
Unfortunately,  a feasible point that satisfies  $T_i \succeq 0$
strictly may not exist.  Hence, if strict feasibility is a concern, one
must use an alternative approach to find reducing certificates in $\approxName$, and potentially
give up on maximum rank certificates.  We propose two alternatives. One, modify the 
strictly feasible variant of $(\star)$ proposed by Lemma 10 of \citet{Lourenco2015} to use $\approxName$. Two, use  a heuristic big-$M$ formulation.  That is, for a large number $M$, solve
\begin{align} \label{sdp:sumOfConesMaxRankSF}
\begin{array}{lll}
 \mbox{\rm maximize}   &  \sum^{|\setName|}_{i=1} \trace{T_i} + M \alpha \  \\
 \mbox{\rm subject to  } & X \in \perpName,\\
&X = \sum^{|\setName|}_{i=1} W_i X_i W^T_i   & \mbox{ \rm  i.e. } X \in \approxName \\
& X_i \succeq T_i  & \forall i\in\{1, \ldots, |\setName|\}\\
&  I \succeq  T_i \succeq (\alpha-1)I &  \forall i\in\{ 1, \ldots, |\setName|\},
\end{array} 
\end{align}
which has a strictly feasible point $(X,X_i,T_i,\alpha) = (0,0,-I,-1)$.  If $M$ is sufficiently large such that $\alpha = 1$ and hence $T_i \succeq 0$ holds  at optimality, then a solution $(X,X_i,T_i,\alpha)$ of \eqref{sdp:sumOfConesMaxRankSF}  yields a solution $(X,X_i,T_i)$ of \eqref{sdp:sumOfConesMaxRank}.

Finally, we mention that strict feasibility of the LP  \eqref{lp:polyMaxRank}
is \emph{not} a concern since exact solutions of LPs can be found efficiently even if strict feasibility fails. That is, finding maximum rank solutions for polyhedral $\approxName$
can be done accurately and efficiently, both in theory and practice, even if strict feasibility of \eqref{lp:polyMaxRank} does not hold.  We illustrate this with examples in Section~\ref{sec:examples}.

\subsection{Explicit modifications of the SDP facial reduction algorithm (Algorithm~\ref{alg:frDual})}
The results of this section are now combined to modify Algorithm~\ref{alg:frDual}. Specifically, we introduce an approximation $\faceName_{i,outer} := U_i \mathcal{C}(\setName_i) U_i^T$ of the face $\faceName_i := U_i \mathbb{S}^{d_i}_{+}U^T_i$
at each iteration $i$, where $\setName_i := \left\{ W_1, \ldots W_{|\setName_i|} \right\}$ is some specified set of rectangular matrices.  A reducing certificate $S_i \in \faceName_{i,outer}^*$
is then found that maximizes the rank of $U_i^T S_i U_i \in \mathcal{C}(\setName_i)^*$  by replacing SDP ($\star$) of Algorithm~\ref{alg:frDual} with the following:
\begin{center}
\fbox{%
        \parbox{0.7\linewidth}{%
\begin{align*}
\begin{array}{lllc}
 \mbox{\rm maximize} &   \sum^{|\setName_i |}_{k=1} \trace{T_k} \  \\
 \mbox{{subject to  }}
 &  U_i U^T_i \cdot S_i > 0  \\
&         U_i^T S_i U_i = \sum^{|\setName_i |}_{k=1} W_k \bar{S}_k W^T_k,
 &\mbox{ i.e., }    U_i^T S_i U_i \in \mathcal{C}(\setName_i)^* \\ 
& \bar{S}_k \succeq T_k  & \forall k\in\{1, \ldots, |\setName_i |\}\\
&  I \succeq  T_k \succeq 0 &  \forall k\in\{ 1, \ldots, |\setName_i |\} \\
& C \cdot S_i  = 0,\;A_j \cdot S_i   = 0 &  \forall j \in \{1,\ldots,m\}.   
\end{array} \; (\star)
\end{align*}
        }%
}
\end{center}
Here, the decision variables are $S_i$ and $T_k,\bar{S}_k$ (for $k\in \{1,\ldots, |\setName_i|\}$) and the first two constraints are equivalent to the condition that $S_i   \in \faceName_{i,outer}^* \setminus \faceName_{i,outer}^{\perp}$
(since $U_i U^T_i \in \relint \faceName_{i,outer})$. To maximize the rank of $U^T_i S_i U_i$, we have applied Corollary~\ref{cor:maxRank}, taking $\perpName$ equal to the subspace $\left\{  U_i^T S U_i : C \cdot S = 0, A_j \cdot S = 0 \right\}$.  
Note the strict inequality  $U_i U^T_i \cdot S_i = \trace  U^T_i S_i U_i > 0$  is satisfied by any non-zero matrix in $\perpName \cap \mathcal{C}(\setName_i)^*$; hence, in practice
one can remove this inequality and instead verify that $U^T_i S_i U_i \ne 0$ holds at optimality, i.e.
one can verify $\perpName \cap \mathcal{C}(\setName_i)^*$ contains a non-zero matrix.

Again note the complexity of solving this modified problem is controlled by $\setName_i$.  When $\setName_i$ contains $d_i \times 1$ rectangular matrices (e.g. $\mathcal{C}(\setName_i)^*$ equals $\D^{d_i}$ or $\DD^{d_i}$), 
the modified problem is a  linear program.  When $\setName_i$ contains $d_i \times 2$ rectangular matrices (e.g. $\mathcal{C}(\setName_i)^*$ equals $\SDD^{d_i}$), it is an SOCP.

\section{Formulation of reduced problems and illustrative examples} \label{sec:simpleEx} 
The facial reduction algorithm for SDP (Algorithm~\ref{alg:frDual}) identifies a face of $\mathbb{S}^n_{+}$
that can be used to  formulate an equivalent SDP. In this section, we   give
simple examples illustrating the basic steps of Algorithm~\ref{alg:frDual} when modified to use approximations  as described in Section~\ref{sec:approach}. We also prove sparsity can be preserved when  certain approximations are used (Proposition~\ref{prop:sparsity}).

\subsection{Formulation of reduced problems} \label{sec:simpProblem}
 Algorithm~\ref{alg:frDual}  identifies a face  $\faceName := U \mathbb{S}^d_{+} U^T$ (where $U \in \mathbb{R}^{n \times d}$ is a fixed
matrix with linearly independent columns and $d \le n$) containing the intersection of $\mathbb{S}^n_{+}$ with an affine subspace $\subspaceName := \{ C - \sum^{m}_{i=1}  y_i A_i : y \in \mathbb{R}^m \}$.
Letting $V \in \mathbb{R}^{n \times n-d}$ denote a matrix whose columns form a basis for $\nullspace U^T$, we can  reformulate the original SDP (reproduced below):
\begin{align*} 
 \begin{array}{lll}
\mbox{\rm maximize} &  b^Ty \\
 \mbox{\rm subject to}  & C - \sum^{m}_{i=1}  y_i A_i   \in \mathbb{S}^n_{+} 
\end{array} 
\end{align*}
explicitly over $\faceName$ as follows:
\begin{align*}
\begin{array}{lll}
\mbox{\rm maximize} &   b^Ty   \\
 \mbox{\rm subject to}  & U^T ( C  - \sum^{m}_{i=1}  y_i A_i ) U   \in \mathbb{S}^{d}_{+} \\ 
			& U^T ( C  - \sum^{m}_{i=1}  y_i A_i ) V = 0 \\
			& V^T ( C  - \sum^{m}_{i=1}  y_i A_i ) V = 0,
\end{array} 
\end{align*}
where we have used  a representation of $\faceName$ given by Lemma~\ref{lem:facesofPSD}.
Here, we see the reduced program is a semidefinite program over $\mathbb{S}^d_{+}$ described by linear equations
and $d \times d$ matrices $U^T C U$ and $U^T A_i U$.

\subsection{Sparsity of reduced problems}
Depending on $U$, the   matrices $U^T C U$ and $U^T A_i U$, though of lower order, may be dense even if $C$ and $A_i$ are sparse.  It is therefore natural to ask how the choice of approximation, discussed in Section~\ref{sec:approach}, affects the structure of $U$ and hence the sparsity of $U^T C U$ and $U^T A_i U$.  It turns out a strong statement can be made when  scaled diagonally-dominant approximations ($\SDD^d$) are used. This statement also applies to diagonal $(\D^d)$ and diagonally-dominant $(\DD^d)$ approximations, since both $\D^d$ and $\DD^d$ are subsets of $\SDD^d$.
 Indeed, a stronger statement can  be made for diagonal approximations. Formally: 
\begin{prop} \label{prop:sparsity}
Let $\faceName$ denote the face identified by the SDP facial reduction algorithm (Algorithm~\ref{alg:frDual}). The following statements hold.
\begin{enumerate}
\item If at each iteration  $U^T_i S_i U_i$ is diagonal, i.e., is in $\D^{d_i}$, then there exists $U  \in \mathbb{R}^{n \times d}$ such that $\faceName= U \mathbb{S}^d_{+} U^T$, where  $U^T X U$ is a principal sub-matrix of $X \in \mathbb{S}^n$ for all $X \in \mathbb{S}^n$. 
\item If at each iteration  $U^T_i S_i U_i$ is scaled diagonally-dominant, i.e., is in $\SDD^{d_i}$, then there exists $U \in \mathbb{R}^{n \times d}$ such that $\faceName= U \mathbb{S}^d_{+} U^T$, where the columns of $U$ have disjoint support. In addition,  $\nnz (U^TX U) \le \nnz (X)$ for all $X \in \mathbb{S}^n$, where $\nnz(\cdot)$ returns the number of non-zero entries of its argument.
\end{enumerate}
\end{prop}
The first statement is trivial to verify.  The next statement is a consequence of Proposition~\ref{prop:kernelsdd}, which implies existence of a matrix $U$ whose columns have disjoint support.  Using disjoint support, the inequality $\nnz ( U^TX U ) \le \nnz (X)$   easily follows.  Also note for diagonal approximations, we can say more about the equations $V^T \subspaceName(y) U=0$ and $V^T \subspaceName(y) V=0$.  In particular, $\nnz(V^TXV) + \nnz(U^TXV) + \nnz(U^TXU) \le \nnz(X)$ for all $X\in \mathbb{S}^n$ since $(U,V)$ can be chosen equal to a permutation matrix.

\subsection{Illustrative Examples}

\subsubsection{Example with diagonal approximations ($\D^d$)}
In this  example, we modify Algorithm~\ref{alg:frDual}  to use diagonal approximations; i.e. at iteration $i$, the face $\faceName_i:= U_i \mathbb{S}^{d_i}_{+} U_i^T$ is approximated by the set $\faceName_{i,outer} = U_i \mathcal{C}(\setName_i) U_i^T$, 
where  $\mathcal{C}(\setName_i)^*$ equals $\D^{d_i}$, the set of $d_i \times d_i$ matrices that are non-negative and diagonal.  A reducing  certificate $S_i$  is found in $\faceName^*_{i,outer}$, the set of matrices $X$ for which $U^T_i X U_i$ is in $\D^{d_i}$.
We apply the algorithm to the following SDP:

\begin{align*}
\begin{array}{lll}
\mbox{\rm Find } & y \in \mathbb{R}^4 & \\
 \mbox{{subject to}} \\
& \subspaceName(y) =  \left(\begin{array}{ccccc} y_1 & 0 & 0 & 0 & 0\\ 
0 & - y_1 & y_2 & 0 & 0\\
0 & y_2 & y_2 - y_3 & 0 & 0\\ 
0 & 0 & 0 & y_3 & 0\\ 
0 & 0 & 0 & 0 & y_4 
\end{array}\right) \in \mathbb{S}^5_{+}.
\end{array}
\end{align*}
Taking $U_0$ equal to the identity matrix and the initial face equal to $\faceName_0 = U_0 \mathbb{S}^5_+ U_0$, we seek a matrix $S_0$ orthogonal to $\subspaceName(y)$ (for all $y$) for  which $U_0^T S_0 U_0$ is  non-negative and diagonal.
An $S_0$  satisfying this constraint and a basis $B$ for $\nullspace U_0^T S_0 U_0$ is given by:
\begin{align*}
S_0  = \left(\begin{array}{ccccc} 
1 & 0 & 0 & 0 & 0\\ 
0& 1 & 0 & 0 & 0\\
 0 & 0 & 0 & 0 & 0\\ 
0 & 0 & 0 & 0 & 0\\ 
0 & 0 & 0 & 0 & 0 \end{array}\right) & \;\;\;\; B = \left(\begin{array}{ccccc} 
 0 & 0 & 0\\ 
 0 & 0 & 0\\
 1 & 0 & 0\\ 
 0 & 1 & 0\\ 
 0 & 0 & 1 \end{array}\right) .
\end{align*}
Taking $U_1 = U_0 B = B$, yields the face $\faceName_1 = U_1 \mathbb{S}_{+}^3 U_1^T$, i.e. the set of PSD matrices in $\mathbb{S}^5_+$ with vanishing first and second rows/cols.

Continuing to the next iteration, we seek a matrix $S_1$ orthogonal to $\subspaceName(y)$ for which $U_1^T S_1 U_1$ is non-negative and diagonal.
An $S_1$  satisfying this constraint and a basis $B$ for $\nullspace U_1^T S_1 U_1$ is given by:
\begin{align*}
S_1  = \left(\begin{array}{ccccc} 
0 & 0 & 0 & 0 & 0\\ 
0 & 0 & -\frac{1}{2} & 0 & 0\\
 0 & -\frac{1}{2} & 1 & 0 & 0\\ 
0 & 0 & 0 & 1 & 0\\ 
0 & 0 & 0 & 0 & 0 \end{array}\right) & \;\;\;\; B = \left(\begin{array}{cccc} 
 0 \\ 
  0 \\ 
 1  \end{array}\right).
\end{align*}
Setting $U_2 = U_1 B$ gives the face $\faceName_2 = U_2 \mathbb{S}^1_{+} U_2^T$, where $U_2 = (0,0,0,0,1)^T$. 

Terminating the algorithm, we now formulate a reduced SDP  over $\faceName_2$. Letting $V$ denote a
basis for $\nullspace U^T_2$ yields:
\begin{align*}
\begin{array}{lll}
\mbox{\rm Find } & y \in \mathbb{R}^4 & \\
 \mbox{\rm subject to  }  & U_2^T \subspaceName(y)   U_2   \in \mathbb{S}^{1}_{+} \\ 
			& U_2^T \subspaceName(y)   V = 0 \\
		        & V^T \subspaceName(y)   V = 0,
\end{array} 
\end{align*}
which simplifies to
\begin{align*}
\begin{array}{lll}
\mbox{\rm Find } & y \in \mathbb{R}^4 & \\
 \mbox{\rm subject to  }  & y_4 \ge 0 \\ 
			& y_1 = y_2 = y_3 = 0.   
\end{array} 
\end{align*}

\paragraph{Existence of reducing certificates}
Lemma~\ref{lem:hp} states that existence of $S_i \in \faceName_{i,outer}^* \setminus \faceName_{i,outer}^{\perp}$ implies   $\subspaceName(y) \cap \relint \faceName_{i,outer}$
is empty. We now verify this fact.  Clearly, $\subspaceName(y)$ is contained in  $\relint \faceName_{0,outer}$
 only if the inequalities 
\[
 y_1 \ge 0 \qquad -y_1 \ge 0
\]
are strictly satisfied, which  cannot hold.  Similarly, $\subspaceName(y)$ is contained in  $\relint \faceName_{1,outer}$ 
only if $y_1 = y_2 = 0$ and the inequalities
\[
 \qquad y_3 \ge 0 \qquad y_2-y_3 \ge 0
\]
are strictly satisfied, which again cannot hold.

\subsubsection{Example with diagonally-dominant approximations ($\DD^{d}$)}
In this  next example, we modify Algorithm~\ref{alg:frDual}  to use diagonally-dominant approximations; i.e. at iteration $i$, the face $\faceName_i:= U_i \mathbb{S}^{d_i}_{+} U_i^T$ is approximated by the set $\faceName_{i,outer} = U_i \mathcal{C}(\setName_i) U_i^T$, 
where  $\mathcal{C}(\setName_i)^*$ equals $\DD^{d_i}$, the set of $d_i \times d_i$ matrices that are diagonally-dominant.  A reducing  certificate $S_i$ is found in $\faceName^*_{i,outer}$, the set of matrices $X$ for which $U^T_i X U_i$ is in $\DD^{d_i}$.
We apply the algorithm to the  SDP

\begin{align*}
\begin{array}{lll}
\mbox{\rm Find } &  y \in \mathbb{R}^3 & \\
 \mbox{subject to  } \\
& \subspaceName(y) =  \left(\begin{array}{cccc}
                  1   &        -y_1    &       0       &    -y_3  \\
               -y_1   &      2y_2 -1   &      y_3      &      0  \\
                   0  &         y_3    &      2 y_1 -1 &   -y_2  \\
               - y_3  &           0    &     -y_2      &      1  
\end{array} \right)  \in \mathbb{S}^4_{+},
\end{array} 
\end{align*}
 and execute just a single iteration of facial reduction. Taking $U_0$ equal to the identity, a matrix $S_0$ orthogonal to $\subspaceName$ for which $U^T_0 S_0 U_0$ is diagonally-dominant and a basis $B$ for $\nullspace U^T_0 S_0 U_0$ is given by
\begin{align*}
S_0 = \left(\begin{array}{cccc} 
1 & 1 & 0 & 0 \\ 
1 & 1 & 0 & 0 \\
0 & 0 & 1 & 1 \\
0 & 0 & 1 & 1 
\end{array}\right)& \;\;  B = \frac{1}{\sqrt{2}} \left(\begin{array}{rr} 
 1 & 0 \\ 
 -1 & 0 \\
 0 & 1 \\ 
 0 & -1 \ 
\end{array}\right).
\end{align*}
Taking $U_1 = U_0 B = B$, yields the face $\faceName_1 = U_1 \mathbb{S}_{+}^2 U_1^T$. (Note the columns of $U$ have disjoint support, a reflection of Proposition~\ref{prop:sparsity}.)

 Terminating the algorithm and constructing the reduced SDP using a matrix $V$ satisfying $\range V = \nullspace U_1^T$ imposes the linear constraints that $y_1 = y_2 = 1, y_3 = 0$; i.e. the reduced SDP has a feasible set consisting of a single point.

\paragraph{Existence of reducing certificates} As was the case in the previous example, existence of a reducing certificate in $\faceName^*_{i,outer} \setminus \faceName_{i,outer}^{\perp}$ implies
emptiness of   $\subspaceName(y) \cap \relint \faceName_{i,outer}$.  We now verify this fact. At the first (and only) iteration, membership of $\subspaceName(y)$ in $\faceName_{0,outer}$ implies  $w_k^T \subspaceName(y) w_k \ge 0$, where $w_k w_k^T$ is any extreme ray of $\DD^4$. 
Taking $w_1 = (1,1,0,0)^T$ and $w_2 = (0,0,1,1)^T$, we have that $\subspaceName(y)$ is contained in $\relint \faceName_{0,outer}$
only if the inequalities
\begin{align*}
\begin{array}{ll}
 w_1^T\subspaceName(y) w_1 =  2 y_2 - 2y_1 \ge 0 \\
w_2^T\subspaceName(y) w_2 =  2 y_1 - 2y_2 \ge 0
\end{array}
\end{align*}
 are  strictly satisfied, which cannot hold.

\section{Recovery of dual solutions} \label{sec:dualRecov}
In this section we address a question that is  relevant to any pre-processing technique based on  facial reduction, i.e. it does not depend in any way on the approximations introduced in Section~\ref{sec:approach}.
Specifically, how (and when) can one recover solutions to the original dual problem? To elaborate, consider the following primal-dual pair\footnote{This designation of primal $(P)$ and dual $(D)$, while standard in
facial reduction literature,  is opposite the convention used by semidefinite solvers such as SeDuMi. We will switch to the convention favored by solvers when we discuss our software implementation in Section~\ref{sec:implementation}. }
for a general conic optimization problem over a closed, convex cone $\coneName$:
\begin{align*} 
\begin{array}[t]{ll} 	
(P): & (D): \\
\;\; \begin{array}{lll}
\mbox{\rm maximize} &  b^Ty  \\
 \mbox{\rm subject to}  & c - Ay   \in \coneName \\ \\
\end{array} \;\;\; & 
\qquad \begin{array}{lll}
\mbox{\rm minimize} &  \langle c, x \rangle\\ 
 \mbox{\rm subject to}  & A^*  x = b &  \\
			& x  \in \coneName^*,
\end{array} 
\end{array}
\end{align*}
and suppose the general facial reduction algorithm (Algorithm~\ref{alg:fr}) is applied to the primal problem $(P)$.
The reduced  primal-dual pair is written over the identified face $\faceName$ and its dual cone $\faceName^*$ as follows:

\begin{align*} 
\begin{array}[t]{ll} 	
(R/P): & (R/D): \\
\;\; \begin{array}{lll}
\mbox{\rm maximize} &  b^Ty \\
 \mbox{\rm subject to}  & c - Ay   \in \faceName \\ \\
\end{array} \;\;\; & 
\qquad \begin{array}{lll}
\mbox{\rm minimize} &  \langle c, x \rangle \\ 
 \mbox{\rm subject to}  & A^*  x = b &  \\
			& x  \in \faceName^*.
\end{array} 
\end{array}
\end{align*}

Since (by construction) $\faceName$  contains $c-Ay$ for any feasible point $y$ of $(P)$,   any solution to $(R/P)$ solves $(P)$.  
On the other hand,   a solution $x$ to $(R/D)$ is not necessarily even a feasible point of  $(D)$
since  $\coneName^* \subseteq \faceName^*$. While  recovering a solution to $(D)$ from a solution to $(R/D)$ may seem in general hopeless, the facial reduction algorithm produces
reducing certificates $s_i \in \faceName^*_i$, where
\[
 \coneName^*  = \faceName^*_0 \subset \faceName^*_1 \subset \cdots \subset \faceName^*_N = \faceName^*,
\]
that can be leveraged to  make recovery possible.  This leads to the following problem statement:

\begin{prob}[Recovery of  dual solutions] \label{prob:prbdsr}
 Given a solution $x$ to $(R/D)$, reducing certificates $s_0,\ldots,s_{N-1}$, i.e. given $s_i$ for which
\begin{align*}
\begin{array}{rcll}
  \langle c,s_i \rangle &=&0  \\  A^*s_i &=& 0    & \\
  s_i &\in& \faceName_i^* \setminus \faceName^{\perp}_i \\
  \faceName_{i+1} &:=& \faceName_{i} \cap s_i^{\perp} \;\;\;\; & ( \mbox{which implies} \;\;\;\;\; \faceName^*_{i+1} = \overline{\faceName^*_i + \lin s_i} )\\
   \faceName_{0} &:=& \coneName, \;\; \faceName := \faceName_N,
\end{array}
\end{align*}
find a solution to $(D)$.
\end{prob}

To solve this problem, we generalize a recovery procedure described in \cite{pataki2016bad} for so-called \emph{well-behaved} SDPs. (See the discussion following \cite[Theorem 5]{pataki2016bad}.) First, we observe that each $s_i$ is a feasible direction for $(R/D)$ that does not increase the dual objective $\langle c, x \rangle$. We
also observe that $\faceName^*_{i+1} = \overline{\faceName^*_i + \lin s_i}$ (since $\coneName$, and hence $\faceName_i$, is closed). This implies if $\faceName^*_i + \lin s_i$ is closed, 
one could, for any $x \in \faceName^*_{i+1}$, find an $\alpha$ such that $x + \alpha s_i$ is in $\faceName^*_i$.  We conclude if $\faceName^*_i + \lin s_i$ were closed for each $i$, then a solution to  $(D)$ could be constructed using a sequence of line searches.   In other words, the following algorithm would successfully recover a solution to $(D)$.

\begin{algorithm}[H] 
 Input: A  solution $x \in \faceName^*$ to the reduced dual  $(R/D)$ and  reducing certificates $s_0 ,\ldots, s_{N-1}$ \\
Output: A solution $x$ to the original dual  $(D)$ or flag indicating failure. \\
    \For{ $i \leftarrow N-1$ \rm down to $0$ }{
\begin{enumerate}
\item Using a line search, find $\alpha$ s.t. $x + \alpha  s_{i} \in \faceName^*_{i}$.
\item If no $\alpha$ exists, return FAIL. Else, set $x \leftarrow x+ \alpha  s_{i}$.
\end{enumerate}
      }
\caption{Recovery of  dual solutions } \label{alg:frPrim_DualRecov}
\end{algorithm} 
\noindent The following properties of this algorithm can be stated immediately:
\begin{lem} \label{lem:success} \
Algorithm~\ref{alg:frPrim_DualRecov} has the following properties:
\begin{enumerate}
\item \emph{Sufficient condition for recovery.} Algorithm~\ref{alg:frPrim_DualRecov} succeeds if $\faceName^*_i + \lin s_i$ is closed for all $i$. 
\item \emph{Necessary condition for recovery.} Suppose $(x,y)$ are optimal solutions to (R/P) and (R/D) with zero duality gap, i.e. $\langle c, x \rangle  =  b^T y$. Then,  Algorithm~\ref{alg:frPrim_DualRecov}
succeeds only if (P) and (D) have solutions with zero duality gap.
\end{enumerate}
\end{lem}
\noindent We note the sufficient condition above \emph{always} holds when $\coneName$ is polyhedral since $\faceName + \mathcal{M}$ is closed  for any subspace $\mathcal{M}$ and face $\faceName$
of a polyhedral cone. On the other hand, $\mathbb{S}^n_{+} + \lin S$ is closed   only if $S$ is zero or positive definite. (To show this, one can use essentially the same argument that proves Lemma~2.2 of \cite{ramana1997strong}.  This lemma shows that for a face $\faceName$ of $\mathbb{S}^n_{+}$, the set $\mathbb{S}^n_{+} + \lin \faceName$ is  closed only if $\faceName = \{0\}$ or $\faceName = \mathbb{S}^n_{+}$.)  Hence, a better sufficient condition for SDP  is desired.

In the next section, we give  a sufficient condition that is also necessary. This condition is specialized to the case $\coneName = \mathbb{S}^n_{+}$ when \emph{one} iteration
of facial reduction is performed.  The restriction to the single iteration case is  imposed so that the condition is easy to state, but it can be extended to the multi-iteration case. We also give a second sufficient condition that is independent of the solution $x$ (Condition~\ref{cond:recEnsure}).

\begin{rem}
 Closure of $\coneName^* +  \lin s$ for $s \in \coneName^*$ has  been studied in other contexts. Borwein and Wolkowicz use 
this condition to simplify their generalized optimality conditions for convex programs (see Remark 6.2 of \cite{borwein1981regularizing}). 
Failure of a related condition, namely closure of $\coneName^* +  \lin \faceName$ for a face $\faceName$ of $\coneName^*$, is used to construct primal-dual pairs with  infinite duality gaps in \cite{tunccel2012strong}.  
\end{rem}

\subsection{A necessary and sufficient  condition for dual recovery} \label{sec:SDPrecov} 

In this section, we give  a necessary and sufficient condition  (Condition~\ref{cond:sdpRecov}) for dual solution recovery that applies when $\coneName = \mathbb{S}^n_{+}$  and a single iteration of facial reduction is performed. 
In this case,
 the primal-dual pair  is given by
\begin{align*} 
\begin{array}{ll} 	
(P-SDP): & (D-SDP): \\
\;\; \begin{array}{lll}
\mbox{\rm maximize} &  b^Ty  \\
 \mbox{\rm subject to}  & C  - \sum^{m}_{i=1}  y_i A_i \succeq 0 
\end{array} &
\qquad \begin{array}{lll} \\
\mbox{\rm minimize} &   C \cdot X  \\ 
 \mbox{\rm subject to}  & A_i \cdot X = b_i  \qquad \forall i \in \{1,\ldots,m\} \\
			& X \succeq 0,
\end{array} 
\end{array}
\end{align*}
where the primal problem is reproduced from Section~\ref{subsec:sdpfr}. The reduced primal-dual pair is over a face $\faceName := \mathbb{S}^n_{+} \cap S^{\perp}$
and its dual cone $\faceName^* =\overline{\mathbb{S}^n_{+} + \lin S}$,
\begin{align*} 
\begin{array}[t]{ll} 	
(R/P-SDP): & (R/D-SDP): \\
\;\; \begin{array}{lll}
\mbox{\rm maximize} &   b^Ty  \\
 \mbox{\rm subject to}  &  \subspaceName(y)  = C  - \sum^{m}_{i=1}  y_i A_i  \\
 &  U^T \subspaceName(y)  U   \in \mathbb{S}^{d}_{+} \\ 
			& U^T \subspaceName(y)  V = 0 \\
			& V^T \subspaceName(y) V = 0 
\end{array} &
\qquad \begin{array}{lll} \\
\mbox{\rm minimize} &   C \cdot X  \\ 
 \mbox{\rm subject to}  & A_i \cdot X = b_i  \qquad \forall i \in \{1,\ldots,m\} \\
			& X = (U,V)  \left( \begin{array}{cc}  W  &  Z  \\  Z^T  & R \end{array} \right) (U,V)^T \\
			& W \in \mathbb{S}^{d}_{+}, R \in \mathbb{S}^{n-d}, Z \in \mathbb{R}^{d \times (n-d)},
\end{array} 
\end{array}
\end{align*}
where $S \in \mathbb{S}^n_{+}$ is a reducing certificate and $(U,V)$ is an invertible matrix satisfying $S=VV^T$ and $\range U = \nullspace S$. Here, the  primal problem is reproduced from Section~\ref{sec:simpProblem}, and the dual problem arises
from a description of $\faceName^*$ given by Lemma~\ref{lem:facesofPSD}.

Algorithm~\ref{alg:frPrim_DualRecov} constructs a solution to  $(D-SDP)$ from a solution $X$ to $(R/D-SDP)$  if and only if $X$ is in $\mathbb{S}^n_{+} + \lin S$.  The following shows this is equivalent to the condition that  $\nullspace W \subseteq \nullspace Z^T$.
We give a   direct proof of this fact, but note it also follows (essentially) by combining \cite[Lemma 3]{pataki2016bad}  with \cite[Lemma 3.2.1]{pataki2000geometry}.  

\begin{lem} \label{lem:dualCone}
Let $(U,V)$ be an invertible matrix for which $\faceName := \mathbb{S}^n_{+} \cap S^{\perp} = U \mathbb{S}^d_{+} U^T$ and  $S = V V^T$. 
 A matrix in the dual cone $\faceName^* = \overline{\mathbb{S}^n_{+} + \lin S}$, i.e., a matrix $X$ of the form
\begin{align} \label{eq:formOfT}
X = (U,V)  \left( \begin{array}{cc}  W  &  Z  \\  Z^T  & R \end{array} \right) (U,V)^T  \qquad \mbox{\rm for some $W \in \mathbb{S}^{d}_{+}, R \in \mathbb{S}^{n-d}, Z \in \mathbb{R}^{d \times n-d}$},
\end{align}
is in $\mathbb{S}^n_{+} + \lin S$ if and only if $\nullspace W \subseteq \nullspace Z^T$.

\begin{proof}
 For the  ``only if'' direction, suppose $X$ is  in $\mathbb{S}^n_{+} + \lin S$, i.e.
 for an $\alpha \in \mathbb{R}$ suppose
\[
X + \alpha V V^T = (U,V)  \left( \begin{array}{cc}  W  &  Z  \\  Z^T  & R+\alpha I \end{array} \right) (U,V)^T \in \mathbb{S}^n_{+}.
\]
Here, membership in $\mathbb{S}^n_{+}$ holds only if $Z^T(I-WW^{\dagger})=0$, where $(I-WW^{\dagger})$ is the orthogonal projector onto $\nullspace W$ (see, e.g. A.5 of \cite{boyd2009convex}). But this implies that
$\nullspace W \subseteq \nullspace Z^T$, as desired.

To see the converse direction, suppose $X$ is such that $Z$ and $W$ satisfy $\nullspace W \subseteq \nullspace Z^T$. The result follows by finding $\alpha$
for which $X+\alpha S \succeq 0$.  We do this by finding an $\alpha_1$ and $\alpha_2$ for which
\begin{align*}
X-VRV^T + \alpha_1 S \succeq 0  \qquad \mbox{and} \qquad VRV^T + \alpha_2S \succeq 0.
\end{align*}
Adding these two inequalities then demonstrates that $X + (\alpha_1+\alpha_2) S \succeq 0$. To find $\alpha_1$, we note that
\[
X-VRV^T + \alpha_1 S = (U,V)  \left( \begin{array}{cc}  W  &  Z  \\  Z^T  & \alpha_1 I \end{array} \right) (U,V)^T.
\]
Taking a Schur complement, the above is PSD if and only if
\[
 W - \frac{1}{\alpha_1} Z  Z^T \succeq 0.
\]
But since $\nullspace W \subseteq \nullspace Z^T$, the matrix $Z  Z^T$ is contained in the face $\mathcal{G} = \left\{ T \in \mathbb{S}^d_{+} : \range T \subseteq \range W \right\}$
where $W$ is in the relative interior of $\mathcal{G}$.  This implies existence of $\alpha_1 > 0 $ for which $W-\frac{1}{\alpha_1} Z  Z^T \in \mathcal{G} \subseteq \mathbb{S}^d_{+}$,
as desired.
To find $\alpha_2$, we note that 
\[
VRV^T+\alpha_2S = V (R+\alpha_2 I)V^T,
\]
where existence of $\alpha_2$ for which $R+\alpha_2 I \succeq 0$ is obvious, completing the proof.

\end{proof}
\end{lem}
\noindent The above characterization of  $\mathbb{S}^n_{+} + \lin S$ yields a necessary and sufficient condition for success of Algorithm~\ref{alg:frPrim_DualRecov} under the assumption that one iteration of facial reduction was performed:
\begin{cond} \label{cond:sdpRecov} 
The solution $X$  to the reduced dual $(R/D-SDP)$ satisfies $\nullspace W \subseteq \nullspace Z^T$.
\end{cond}


\noindent The following example illustrates success and failure of Condition~\ref{cond:sdpRecov}.

\begin{ex}  \label{ex:exampleRecov}
Consider the following primal-dual pair:
\begin{align*}
\begin{array}{ll} 
\begin{array}{ll}
         \mbox{\rm maximize } &  y_3 + 2y_2  \\ \mbox{\rm subject to  } \\
   	   & \subspaceName(y) =  \left(\begin{array}{ccc} {y_1} & y_2 & 0\\ 
				        y_2 & -{y_3} & {y_2}\\ 0 & {y_2} & {y_3}  \end{array} \right)  \succeq 0  \\
	\end{array} 
\end{array} & \;\;\;\;
      \begin{array}{lll}
         \mbox{\rm minimize }    & 0
	   \\ \mbox{\rm subject to } 
&	  x_{33} - x_{22} = -1 \\
&	  x_{12} + x_{21} + x_{23} +x_{32} = -2 \\
&	  x_{11} = 0 \\
	& X \succeq 0  \\
	\end{array} 
\end{align*}
and let $S = VV^T$, with $V=(e_2,e_3)$.  Clearly, $S$ is a reducing certificate defining a face $\faceName := \mathbb{S}_{+}^n \cap S^{\perp} = U\mathbb{S}^{1}_{+} U^T$ for $U = e_1 = (1,0,0)^T$.  
Rewriting the primal-dual pair over $\faceName$ and $\faceName^*$ gives:
\begin{align*}
\begin{array}{ll}
\begin{array}{ll}
         \mbox{\rm maximize } &  y_3 + 2y_2  \\ \mbox{\rm subject to} \\
         &  V^T \subspaceName(y) V = 0 \\
         &  U^T \subspaceName(y) V = 0 \\
         & U^T \subspaceName(y) U \succeq 0 \\
   	 &  \subspaceName(y) =   \left(\begin{array}{ccc} {y_1} & y_2 & 0\\ 
				        y_2 & -{y_3} & {y_2}\\ 0 & {y_2} & {y_3}  \end{array} \right) 
	\end{array} 
\end{array}  & 
      \begin{array}{lll}
         \mbox{\rm minimize }    & 0
	   \\ \mbox{\rm subject to} 
&	  x_{33} - x_{22} = -1 \\
&	  x_{12} + x_{21} + x_{23} +x_{32} = -2 \\
&	  x_{11} = 0 \\
&  X \in \mathbb{S}^3,  U^T X U  = x_{11} \ge 0 \\ \\ \\ \\
\end{array}
\end{align*}
A solution to the dual problem that satisfies Condition~\ref{cond:sdpRecov} is given by:
\begin{align*}
X = \left(\begin{array}{ccc} 0  & 0 & 0\\ 0 & 0 & -1 \\ 0 & -1 & -1 \end{array}\right).
\end{align*}
To see the condition is satisfied, note $Z = (x_{12},x_{13}) = (0,0)$ and $W = x_{11} = 0$.  Hence, $\nullspace Z^T$  contains (indeed, equals) $\nullspace W$.
We therefore see that solution recovery succeeds, i.e. for (say) $\alpha =2$:
\[
 X + \alpha S = \left( \begin{array}{ccc} 0 & 0 & 0 \\ 0 & \alpha & -1 \\ 0 & -1 & \alpha-1 \end{array} \right) \succeq 0.
\]
 Conversely, the following solution fails Condition~\ref{cond:sdpRecov}:
\begin{align} \label{sol:fail}
X =  \left(\begin{array}{ccc} 0  & -1 & 0\\ -1 & 0 & 0 \\ 0 & 0 & -1 \end{array}\right).
\end{align}
Here, $Z = (-1,0)$ and $W = 0$.  Hence, $\nullspace Z^T = \{0\}$   does not contain $\nullspace W = \mathbb{R}$ and recovery must fail.
In other words, there is no $\alpha$ for which 
\[ 
X + \alpha S = \left( \begin{array}{ccc} 0 & -1 & 0 \\ -1 & \alpha & 0 \\ 0 & 0 & \alpha-1 \end{array} \right) \succeq 0,
\]
which is easily seen.
\end{ex}

\subsubsection{Strong duality is not sufficient for dual recovery}
An additional observation can be made about Example~\ref{ex:exampleRecov}.  As we observed in Lemma~\ref{lem:success}, zero duality gap between the original primal-dual pair $(P)$ and $(D)$
is a necessary condition for recovery to succeed when the reduced primal-dual pair $(R/D)$ and $(R/D)$ has zero duality gap.  Example~\ref{ex:exampleRecov} shows
this is \emph{not} a sufficient condition when $\coneName = \mathbb{S}^n_{+}$. Here, both the original  and the reduced primal-dual pairs have zero duality gap yet successful recovery depends on
 the specific solution found for the reduced dual $(R/D-SDP)$.
This is summarized below:
\begin{cor}
  The dual solution recovery procedure of Algorithm~\ref{alg:frPrim_DualRecov} can fail even  if both the original primal-dual pair $(P)$ and $(D)$ 
and the reduced primal-dual pair $(R/P)$ and $(R/D)$  have zero duality gap.
\end{cor}

\subsubsection{Ensuring successful dual recovery} \label{sec:ensureRec}
Condition~\ref{cond:sdpRecov} lets one determine if recovery is possible by a simple null space computation.  Unfortunately, this check
must be done \emph{after} the dual problem $(R/D-SDP)$ has been solved as it depends on the specific  solution that is obtained.   In this section, we give a simple sufficient 
condition that can be checked \emph{prior} to solving $(R/D-SDP)$.  If this condition is satisfied, a modification of the reduced primal-dual pair 
can be performed to guarantee successful recovery, independent of the solution found for   $(R/D-SDP)$;
the explicit modification is given by $(R/P-SDP-2)$ and $(R/D-SDP-2)$.

The idea is simple: when can one assume $Z = 0$  and hence ensure Condition~\ref{cond:sdpRecov} holds  without loss of generality?  By \cite[Theorem 5]{pataki2016bad},  this assumption can be made if   $(R/P-SDP)$ satisfies Slater's condition and $(P-SDP)$ is \emph{well-behaved}---where
$(P - SDP)$ is well-behaved  if, for \emph{all} cost vectors $b$, the SDPs  $(P - SDP)$ and  $(D - SDP)$ have no duality gap and $(D - SDP)$  attains it optimal value when it is finite.    It turns out we can assume $Z=0$ under a related but purely linear-algebraic condition  inspired by a  characterization of well-behaved SDPs \cite[Theorem 3]{pataki2016bad}.    The condition and  statement follow.
\begin{cond} \label{cond:recEnsure}
The equations of $(R/P-SDP)$ have the following property:  
\begin{align*} 
  \left\{ y \in \mathbb{R}^m : V^T \subspaceName(y) V = 0 \right\} &= \left\{ y \in \mathbb{R}^m : V^T \subspaceName(y)  V = 0, V^T \subspaceName(y) U = 0 \right\},
\end{align*}
that is, $V^T \subspaceName(y)  V = 0$ implies  $V^T \subspaceName(y)  U = 0$.
\end{cond}
\begin{prop} \label{prop:blockdiag}
Suppose Condition~\ref{cond:recEnsure} holds.  If $(R/D-SDP)$ has  an optimal solution, then it has an optimal solution  with $Z = 0$.

\begin{proof}
Let $X$ be an optimal solution to $(R/D-SDP)$, which, for some $W \in \mathbb{S}^d_{+}$, $Z \in \mathbb{R}^{d \times (n-d)}$, and $R \in \mathbb{S}^{n-d}$ satisfies 
\begin{align} \label{eq:Xdecomp}
X = (U,V)  \left( \begin{array}{cc}  W  &  Z  \\  Z^T  & R \end{array} \right) (U,V)^T.
\end{align}
 We will   construct  a new solution $\hat X$ by setting $Z$ to zero and replacing $R$ with $R+\hat R$ for a particular $\hat R$.

Towards this, we first show existence of $X$ implies the set $ \left\{ y \in \mathbb{R}^m : V^T \subspaceName(y) V = 0 \right\}$ is non-empty. If it were empty, then, by Farka's lemma, there exists $\tilde R$ satisfying $\tilde R \cdot (V^TA_i V) = 0$ and $\tilde R \cdot (V^T C  V) < 0$, which  implies
\[
\tilde X = X+ (U,V)  \left( \begin{array}{cc}  0  &  0  \\  0  & \tilde R \end{array} \right) (U,V)^T 
\]
 is a feasible point of $(R/D-SDP)$ with strictly better cost, contradicting optimality of $X$. Hence, there exists $y_0 \in \left\{ y \in \mathbb{R}^m : V^T \subspaceName(y) V = 0 \right\}$.

Now, consider the linear maps $L_1  : \mathbb{R}^m \rightarrow \mathbb{S}^{ n-d} $ and $L_2  : \mathbb{R}^m \rightarrow \mathbb{R}^{d \times (n-d)} $ and
corresponding adjoint maps $L^*_1  : \mathbb{S}^{ n-d} \rightarrow \mathbb{R}^m $ and $L^*_2  : \mathbb{R}^{d \times (n-d)}   \rightarrow \mathbb{R}^m$
defined via
\begin{align*}
L_1(y) &=  \sum^m_{i=1} y_i (V^T A_i  V), \;\; &L_1^*(R)  &= \left( (V^T A_1 V) \cdot  R, \cdots, (V^T A_m V) \cdot R \right)^T  & \\
L_2(y)  &=  \sum^m_{i=1} y_i (U^T A_i V)   \;\; &L_2^*(Z)  &= \left( (U^T A_1 V) \cdot  Z  , \cdots, (U^T A_m V) \cdot Z \right)^T,
\end{align*}
where $\mathbb{R}^{d \times (n-d)}$ and $\mathbb{S}^{ n-d} $ are equipped with trace inner-product  $P \cdot Q := \trace P^TQ$. With these definitions,   $L_1(y) = V^TCV$ iff $V^T \subspaceName(y)  V = 0$
and $L_2(y) = U^TCV$ iff $V^T \subspaceName(y)  U = 0$. Further, 
$X$ of form \eqref{eq:Xdecomp} satisfies the equations $A_i \cdot X = b_i$ for $(R/D-SDP)$ iff
\[
 L_1^*(R)+ 2L_2^*(Z) = b - \left( (U^TA_1U ) \cdot W, \ldots, (U^TA_mU) \cdot W \right)^T .
\]
Now suppose   Condition~\ref{cond:recEnsure} holds. Given existence of $y_0$, it follows that $\nullspace L_1 \subseteq \nullspace L_2$---otherwise,  we could construct solutions to $L_1(y) = V^TCV$ that do not solve  $L_2(y)=U^T C V$, a contradiction of Condition~\ref{cond:recEnsure}.  But $\nullspace L_1 \subseteq \nullspace L_2$
holds if and only if    $\range L^*_1 \supseteq \range L^*_2$.
Hence, we can find a $\hat R$ satisfying
\begin{align} \label{eq:range}
L^*_1(\hat R) = 2 L^*_2(Z),
\end{align}
which implies the matrix
\[
\hat X = (U,V)  \left( \begin{array}{cc}  W  &  0  \\  0  & R+\hat R \end{array} \right)  (U,V)^T 
\]
satisfies $A_i \cdot \hat X = b_i$. Since $W \in \mathbb{S}^d_{+}$, it follows $\hat X$ is feasible for $(R/D-SDP)$.  

We now show $C\cdot X = C \cdot \hat X$, proving $\hat X$ is also optimal. For this, it suffices to show $V \hat R V^T \cdot C = 2 U Z V^T \cdot C$.  Since  $L_1(y_0) = V^TCV$ and, by Condition~\ref{cond:recEnsure}, $L_2(y_0)=U^T C V$, we conclude
\begin{align*}
\langle L^*_1(\hat R), y_0 \rangle &= \langle \hat R, L_1(y_0) \rangle  =  \langle \hat R, V^TCV  \rangle = V \hat R V^T \cdot C, \\
\langle  L^*_2( Z), y_0 \rangle &= \langle  Z, L_2(y_0) \rangle  =  \langle Z, U^TCV  \rangle = U Z V^T \cdot C,
\end{align*}
which, by \eqref{eq:range}, shows $V \hat R V^T \cdot C = 2 U Z V^T \cdot C$, as desired.
\end{proof}
\end{prop}

We conclude one can fix $Z$ to zero in $(R/D-SDP)$   and omit  the equations  $V^T \subspaceName(y) U = 0$ from $(R/P-SDP)$ under Condition~\ref{cond:recEnsure}.   This leads to a modified primal-dual pair:
\begin{align*} 
\begin{array}[t]{ll} 	
(R/P-SDP-2): & (R/D-SDP-2): \\
\;\; \begin{array}{lll}
\mbox{\rm maximize} &  b^Ty  \\
 \mbox{\rm subject to}  &  \subspaceName(y)  = C  - \sum^{m}_{i=1}  y_i A_i  \\
 &  U^T \subspaceName(y)  U   \in \mathbb{S}^{d}_{+} \\ 
			& V^T \subspaceName(y) V = 0 
\end{array} &
\qquad \begin{array}{lll} \\
\mbox{\rm minimize} &   C \cdot X  \\ 
 \mbox{\rm subject to}  & A_i \cdot X = b_i  \qquad \forall i \in \{1,\ldots,m\} \\
			& X = (U,V)  \left( \begin{array}{cc}  W  &  0  \\  0  & R \end{array} \right) (U,V)^T \\
			& W \in \mathbb{S}^{d}_{+}, R \in \mathbb{S}^{n-d},
\end{array}
\end{array}
\end{align*}
where any solution to $(R/P-SDP-2)$ solves the original primal  $(P-SDP)$ and any solution to $(R/D-SDP-2)$ satisfies  Condition~\ref{cond:sdpRecov} (by construction).  Note given  a solution of $(R/D-SDP-2)$, Algorithm~\ref{alg:frPrim_DualRecov} recovers a solution to $(D-SDP)$  that is block-diagonal (i.e., satisfies $Z=0$); identical block-diagonal structure was established for well-behaved SDPs by \cite[Theorem 5]{pataki2016bad}.

\paragraph{Comparison with well-behavedness}  
 
We now illustrate differences between Condition~\ref{cond:recEnsure}  and well-behavedness of $(P-SDP)$. Suppose $(R/P-SDP)$ is constructed using one iteration of facial reduction.  From \cite[Theorem 3]{pataki2016bad}, it follows  Condition~\ref{cond:recEnsure}  and well-behavedness of $(P-SDP)$ are equivalent if $(R/P-SDP)$  satisfies Slater's condition. The following examples show this equivalence can fail if Slater's condition does not hold.

\begin{ex}[A well-behaved SDP and failure of Condition~\ref{cond:recEnsure}]
 Consider the following SDP:
 \begin{align*}
 \begin{array}{llc}
 \mbox{\rm maximize } & b^Ty  \\
  \mbox{\rm{subject to }} \\
 & \subspaceName(y) =  \left(\begin{array}{ccccc} y_1 & 0 & y_2 & 0 \\ 
 0 & - y_1 & 0 & 0 \\
 y_2 & 0 & y_2 & 0 \\ 
 0 & 0 & 0 & - y_2 
 \end{array}\right) \in \mathbb{S}^4_{+}.
 \end{array}
 \end{align*}
 The feasible set is given by $y_1 = y_2 = 0$, and a dual optimal solution  by a non-negative diagonal matrix $X$ satisfying  $b_1 = x_{22}-x_{11}$ and $b_2 = x_{44}-x_{33}$, which clearly exists for all $b$.  Hence, this SDP is well-behaved. 
 
  For the rank two reducing certificate $S =  e_1e_1 + e_2 e_2^T$, we have that $(R/P -SDP)$ takes the form:
  \begin{align*}
 \begin{array}{llcc}
 \mbox{\rm maximize } & b^Ty & \\
  \mbox{\rm{subject to }} \\
 & U^T \subspaceName(y)U =  \left(\begin{array}{ccccc} y_2 & 0 \\ 
 0 & - y_2 
 \end{array}\right)  \in \mathbb{S}^2_{+}, 
 & V^T \subspaceName(y)U =  \left(\begin{array}{ccccc} y_2 & 0 \\ 
 0 & 0
 \end{array}\right)  = 0_{2 \times 2}, \\ 
 &  V^T \subspaceName(y)V =  \left(\begin{array}{ccccc} y_1 & 0 \\ 
 0 & -y_1
 \end{array}\right)  = 0_{2 \times 2} 
 \end{array}
 \end{align*}
 where $V=(e_1,e_2)$, $U=(e_3,e_4)$ and $S=VV^T$.  Clearly $U^T \subspaceName(y)U \succeq 0$ cannot be strictly satisfied, hence $(R/P-SDP)$ fails Slater's condition. In addition,    Condition~\ref{cond:recEnsure}  fails, i.e.,
     \[
    \left\{ y \in \mathbb{R}^m : V^T \subspaceName(y) V = 0 \right\} \ne \left\{ y \in \mathbb{R}^m : V^T \subspaceName(y)  V = 0, V^T \subspaceName(y) U = 0 \right\}.
   \] 
 Note failure of Slater's condition occurs because we did not use the full rank reducing certificate $S=I$.
\end{ex} 
 
 \begin{ex}[An SDP not well-behaved and success of Condition~\ref{cond:recEnsure}]
 The following SDP, based on Example 1 of  \citet{pataki2016bad},
 has a dual optimal value that is unattained when $b=(1,0,0)^T$ and is hence not well-behaved:
 \begin{align*}
 \begin{array}{llc}
 \mbox{\rm maximize } & b^Ty  \\
  \mbox{\rm{subject to }} \\
 & \subspaceName(y) =  \left(\begin{array}{ccccc} 1 & -y_1 & 0 & 0 \\ 
 -y_1 & y_3 & 0 & 0 \\
 0 & 0 & y_2 & y_3 \\ 
 0 & 0 & y_3 & - y_2 
 \end{array}\right) \in \mathbb{S}^4_{+}.
 \end{array}
 \end{align*}
 For the rank two reducing certificate $S =  e_3e_3 + e_4 e_4^T$, we have that  $(R/P-SDP)$ takes the form:
  \begin{align*}
 \begin{array}{llcc}
 \mbox{\rm maximize } & b^Ty  \\
  \mbox{\rm{subject to}} \\
 & U^T \subspaceName(y)U =  \left(\begin{array}{ccccc} 1 & -y_1 \\ 
 -y_1 & y_3 
 \end{array}\right)  \in \mathbb{S}^2_{+},
 & V^T \subspaceName(y)U =  \left(\begin{array}{ccccc} 0 & 0 \\ 
 0 & 0
 \end{array}\right)  = 0_{2 \times 2}, \\ 
 &  V^T \subspaceName(y)V =  \left(\begin{array}{ccccc} y_2 & y_3 \\ 
 y_3 & -y_2
 \end{array}\right)  = 0_{2 \times 2},
 \end{array}
 \end{align*}
 where $V=(e_3,e_4)$ and $U=(e_1,e_2)$. Since $y_3 = 0$ if $y$ is  feasible, Slater's condition fails.  Condition~\ref{cond:recEnsure}, on the other hand, holds given that $V^T \subspaceName(y)U = 0$ imposes no constraints on $y$. This is despite the fact the SDP is not well-behaved. Also note  $(R/P-SDP)$  fails Slater's condition even though $S$ is a reducing certificate of maximum rank. 
\end{ex}

\subsection{Recovering solutions to an extended  dual}
We close this section by discussing recovery for an alternative dual program intimately related to facial reduction---a so-called  \emph{extended dual} \cite{pataki2013simple}.  
For SDP, this dual is a slight variant of the Ramana dual \cite{ramana1997exact},
which was related to facial reduction in \cite{ramana1997strong}.  

A solution to an extended dual  carries the same information as a solution  to the reduced dual $(R/D)$ and a sequence of reducing certificates
used to identify a face. However, such a solution  allows one to certify optimality of the primal problem $(P)$ without retracing 
the steps of the facial reduction algorithm (to verify validity of each reducing certificate)---one simply checks that a solution to an extended dual and a candidate solution to $(P)$ have zero duality gap.   

Extended duals can be defined for cones $\coneName$ that are \emph{nice}  \cite{pataki2013simple}, but we will limit  discussion to the case when $\coneName = \mathbb{S}^n_{+}$. The extended
dual considered is based on three  key facts. 
\begin{lem} The following statements are true:
\begin{enumerate}
\item For any face $\faceName$ of $\mathbb{S}^n_{+}$, $\faceName^* = \mathbb{S}^n_{+} + \faceName^{\perp}$. 
\item If $\faceName = \mathbb{S}^n_{+} \cap S^{\perp}$ for $S \in \mathbb{S}^n_{+}$, then 
\[
\faceName^{\perp} = \left\{ W + W^T : \left(\begin{array}{cc} S & W \\ W^T & \alpha I \end{array}\right) \succeq 0  \mbox{ for some } \alpha \in \mathbb{R}  \right\}.
\]
\item Let $\faceName_0 := \mathbb{S}^n_{+}$ and consider the chain of faces defined by matrices $S_i$
\[
 \faceName_{i+1} :=  \faceName_i  \cap S^{\perp}_i, 
\]
where $S_i$ is in $\faceName^*_i$, i.e. $S_i = \bar{S}_i + V_i$ for $\bar{S}_i \in \mathbb{S}^n_{+}$ and $V_i \in \faceName^{\perp}_i$.  The following relationship holds:
\[
 \faceName_{i+1} =\mathbb{S}^n_{+} \cap ( \sum^{i}_{j=0} \bar{S}_j)^{\perp}.
\]
\begin{proof}
The first statement holds because $\mathbb{S}^n_{+}$ is a  \emph{nice} cone \cite{pataki2013simple}.  The other statements are shown by Proposition 1 and Theorem 3 of \cite{pataki2013simple}.
\end{proof}

\end{enumerate}
\end{lem}
\noindent  Using these facts, the extended dual  considered simultaneously identifies a chain of faces $\faceName_1,\ldots,\faceName_N$  (where $N$ can be chosen to equal the length of the longest chain of faces of $\mathbb{S}^n_{+}$) and
a solution $X \in \faceName_N^*$  to the reduced dual $(R/D-SDP)$ formulated over $\faceName_N^*$. 
It is given below as an optimization problem over $X,\bar{X},W_N, S_i, \bar{S}_i,W_i, \alpha_i$:
\begin{align*}
\begin{array}{l} (EXT/D-SDP) \\
      \begin{array}{rrl}
         \mbox{ minimize }    & C \cdot X
	   \\ \mbox{ subject to} 
	& A_j \cdot X = b_j   \\
       &  C \cdot S_i = 0, \; A_j \cdot S_i = 0 & (\mbox{i.e. $S_i^{\perp}$ contains $\subspaceName$})  \\
        & X = \bar{X} + W_N + W^T_N & (\mbox{i.e. $X \in \mathbb{S}^n_{+} + \faceName_N^{\perp} = \faceName^{*}_N$}) \\
	 & S_i = \bar{S_i} + W_{i} + W_{i}^T  & (\mbox{i.e. $S_i \in \mathbb{S}^n_{+} + \faceName_{i}^{\perp} = \faceName^{*}_{i}$}) \\  
	 & \left(\begin{array}{cc} \sum^i_{j=0} \bar{S}_j & W_{i+1} \\ W_{i+1}^T & \alpha_i I \end{array} \right) \succeq 0 & (\mbox{i.e. $W_{i+1} +W_{i+1}^T \in \faceName_{i+1}^{\perp}$}) \\
             & \bar{S}_i \succeq 0, \bar{X} \succeq 0, W_0 = 0,
\end{array}
\end{array}
\end{align*}
where $i$ ranges from $0$ to $N-1$ and $j$ ranges from $1$ to $m$ (indexing $m$   linear equations  $A_j \cdot X = b_j$).

\paragraph{Recovering a solution}
Suppose $\faceName_i = U_i \mathbb{S}_{+}^{d_i} U_i^T$ for $i=0,\ldots,N$ is a sequence of  faces identified by an SDP facial reduction procedure (e.g. Algorithm~\ref{alg:frDual}, with or without the modifications of Section \ref{sec:approach}) 
suitably padded so that the length of the sequence is $N$, i.e. $\faceName_0,\ldots,\faceName_M = \mathbb{S}^n_{+}$
for some $M < N$.  Let  $S_i  \in \faceName^*_i$ be the corresponding sequence of reducing certificates (similarly padded with zeros) and let $X$ be a solution to $(R/D-SDP)$. One can construct a feasible point to $(EXT/D-SDP)$ by  decomposing $S_i$ (and similarly $X$) into the form  $S_i = \bar{S}_i + W_i + W^T_i$, for $\bar{S}_i \in \mathbb{S}^n_{+}$ and $W_i + W^T_i \in \faceName_{i}^{\perp}$. Supposing 
$U_i$ has orthonormal columns,  this can be done by taking:
\begin{align*}
\begin{array}{lll}
 \bar{S}_i =  U_i U^T_i S_i U_i U^T_i & W_i = \frac{1}{2} (  S - \bar{S}_i) & \forall i \in \{0,\ldots,N-1\} \\
  \bar{X} = U_N U^T_N X U_N U_N^T & W_N = \frac{1}{2} (  X - \bar{X}) .
\end{array}
\end{align*}
One can then pick $\alpha_i$ (individually) until the relevant semidefinite constraint is  satisfied. The feasible point produced by this procedure is optimal
if the reduced primal-dual pair over $\faceName_N$ and $\faceName_N^*$ has no duality gap. This of course occurs if the reduced primal problem over $\faceName_N$ is 
strictly feasible (i.e. the unmodified version of Algorithm~\ref{alg:frDual} is run to completion).

\section{Implementation} \label{sec:implementation}

\begin{figure}
\centering
\includegraphics[scale=.4]{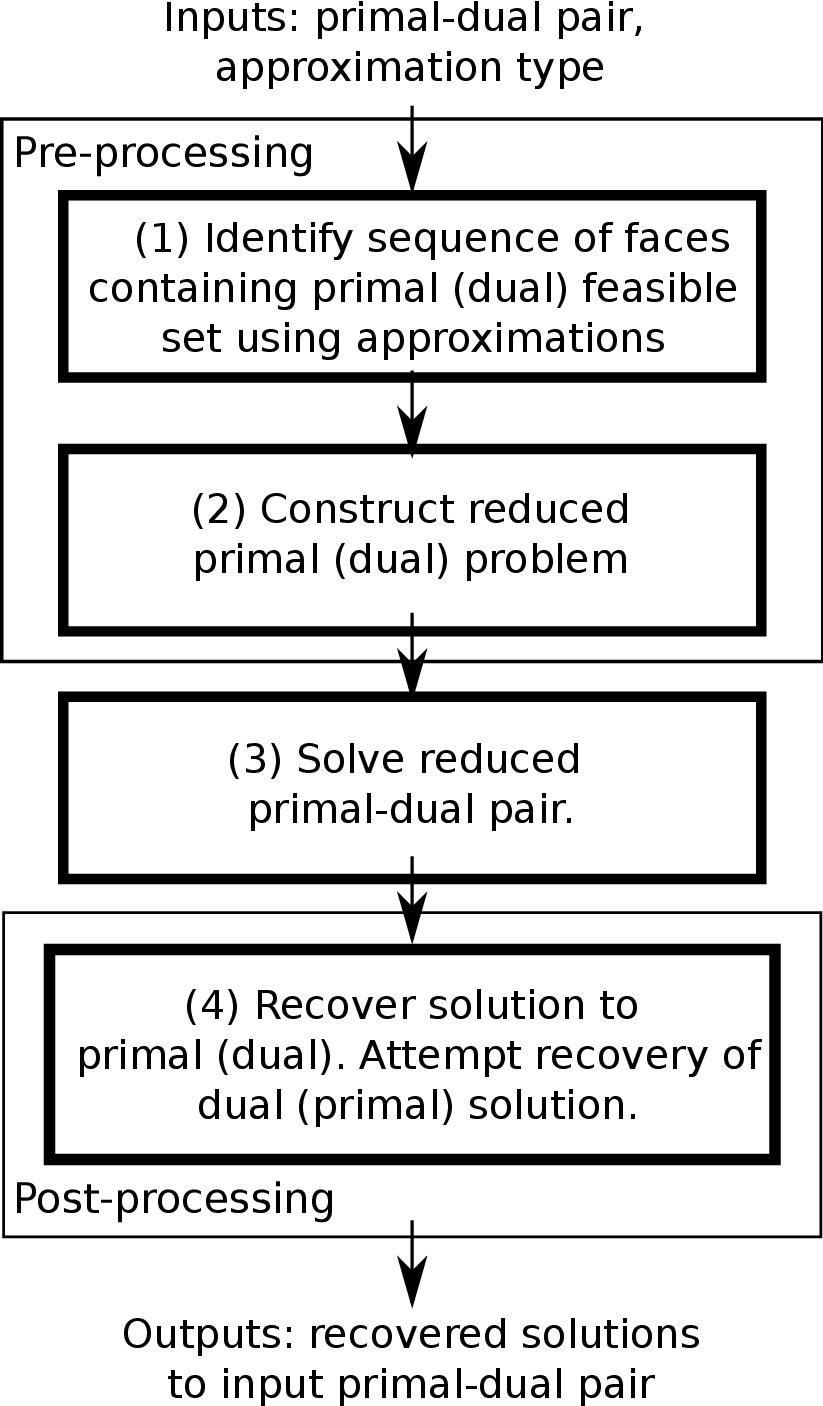}
\caption{Flow  of  MATLAB implementation}
\label{fig:flow}
\end{figure}

The discussed techniques have been implemented as a suite of MATLAB scripts we dub {\tt frlib},  available at
at \tt{www.mit.edu/\textasciitilde fperment}\rm.   The basic work flow is depicted in Figure~\ref{fig:flow}. The implemented code takes as input a  primal-dual SDP \emph{pair}  
and can reduce  (using  suitable variants of Algorithm~\ref{alg:frDual})
either the primal problem or the dual. This is an important feature  since either the primal or the dual  may model the problem of interest.

\subsection{Input formats}
The implementation takes in SeDuMi-formatted inputs  {\tt A,b,c,K}, where {\tt A,b,c}, define
the subspace constraint and objective function and {\tt K} specifies the sizes of the semidefinite constraints   \cite{sturm1999using}. 
Conventionally, the primal problem  described by {\tt A,b,c,K} refers to an SDP  defined by equations $A_i \cdot X = b_i$. 
Similarly, the dual problem described by {\tt A,b,c,K}  refers to an SDP defined by generators $C-\sum_i y_i A_i$.  
While our implementation and the following discussion follow this convention,   the  \emph{opposite}  convention was used in previous sections (e.g.  $(SDP-P)$  and  $(SDP-D)$ in   Section~\ref{sec:SDPrecov}).

\subsection{Reduction of the primal problem} \label{sec:primRed}
Given {\tt A,b,c,K}; the following syntax is used to reduce the primal problem, solve the reduced primal-dual pair, and recover solutions to the original primal-dual pair
via our implementation:
\begin{center}
\begin{verbatim}
prg = frlibPrg(A,b,c,K);
prgR = prg.ReducePrimal(`d');
[x_reduced,y_reduced] = sedumi(prgR.A, prgR.b, prgR.c, prgR.K);
[x,y,dual_recov_success] = prgR.Recover(x_reduced,y_reduced);
\end{verbatim}
\end{center}
The call to {\tt prg.ReducePrimal} reduces the primal problem using diagonal ({ \tt `d' }) approximations by executing a variant of Algorithm~\ref{alg:fr}.
To find reducing certificates, it solves a series of LPs (defined by the diagonal approximation) that
can be solved using a handful of supported solvers. The returned object {\tt prgR} has member variables
\begin{center}
{\tt prgR.A, prgR.b, prgR.c, prgR.K},
\end{center}
which describe the reduced primal-dual pair. For a single semidefinite constraint, this reduced primal-dual pair
is given by:
\begin{align*}
 \begin{array}{lll}
\mbox{\rm minimize} &  C \cdot U \hat{X}U^T  \\ 
 \mbox{\rm subject to}  & A_i \cdot U \hat{X}U^T  = b_i \;\; \forall i \in \{1,\ldots,m\} \\
			& \hat{X}  \in \mathbb{S}^d_{+}
\end{array} \qquad
\begin{array}{lll}
\mbox{\rm maximize} &  b^Ty  \\
 \mbox{\rm subject to}  & U^T (C - \sum^{m}_{i=1}  y_i A_i) U  \in \mathbb{S}^d_{+}, \\ \\
\end{array} 
\end{align*}
where $U \mathbb{S}_{+}^d U^T$ is a face identified by {\tt prg.ReducePrimal}. 
The reduced primal and its dual are solved by calling SeDuMi.

The primal solution {\tt x\_reduced} returned by  SeDuMi represents  an optimal $\hat{X}$. 
The function {\tt prgR.Recover} computes from $\hat{X}$ a solution $U \hat{X}U^T$ to the original primal problem.  It then attempts to 
find a solution to the original dual using a variant of the recovery procedure
described in Section~\ref{sec:dualRecov}  (Algorithm~\ref{alg:frPrim_DualRecov}). The flag {\tt dual\_recov\_success}  indicates  success of this recovery procedure.

\subsection{Reduction of the dual problem}
The above syntax can be modified to reduce the dual problem described by {\tt A,b,c,K}. This is done replacing the relevant line above with:
\begin{verbatim}
prgR = prg.ReduceDual(`d');
\end{verbatim}
As above, the object {\tt prgR}  contains a description of the primal-dual pair which, for a single semidefinite constraint,
is given by the SDPs $(R/P-SDP)$ and $(R/D-SDP)$ in Section~\ref{sec:SDPrecov} (where, recalling our earlier convention, the label of ``primal'' and ``dual'' is reversed). 
With {\tt prgR} created in this manner, a call to {\tt prgR.Recover} (though syntactically identical) now returns a solution to the original \emph{dual}
and attempts to recover a solution to the original \emph{primal} using Algorithm~\ref{alg:frPrim_DualRecov}. In other words, a  call of the form
\begin{verbatim}
[x,y,prim_recov_success] = prgR.Recover(x_reduced,y_reduced);
\end{verbatim}
returns a solution {\tt y}  to the original dual problem  and attempts to recover  a solution {\tt x}  to the original  primal problem.
The flag {\tt prim\_recov\_success} indicates successful recovery of {\tt x}.

\subsection{Solution recovery}
As suggested by the flags  {\tt prim\_recov\_success}  and {\tt dual\_recov\_success}  in the preceding examples, solution recovery is only guaranteed for the  problem that is reduced, i.e.
if the primal (resp. dual) is reduced, recovery of the original dual (resp. primal) may fail for reasons discussed  in Section~\ref{sec:SDPrecov}. Thus, it is important to reduce the primal only if it is the problem of interest,
and similarly for the dual.

\section{Examples} \label{sec:examples}
This section gives larger examples that illustrate  effectiveness of our method. For each example, the same type of approximation (e.g. diagonal or diagonally-dominant) is  used at each facial reduction iteration. Many  examples are also over products of  cones, e.g. $\coneName = \mathbb{S}^{n_1} \times \mathbb{S}^{n_2} \times \cdots \times \mathbb{S}^{n_k}$.
In these cases, we use the same type of approximation for each cone $\mathbb{S}^{n_i}$. 
For some examples, the feasible set is defined by equations $A_i \cdot X = b_i$; for these examples, reduced SDPs are formulated  using  obvious variants of our method.  Also, for a reduced SDP of the form $(R/P-SDP)$, the equations  $V^T \subspaceName(y) U = 0$ and $V^T \subspaceName(y) V = 0$ are eliminated before solving the SDP. 
For each example, we report one or more of the following items $(1-4)$:

\paragraph{1) Complexity parameters and sparsity}
For each example, we report a list of numbers describing the size and sparsity of the problem, denoted
\[
n;r;\nnz.
\]
Here, $n$ gives the size(s) of the psd cone(s)  and $r$ the dimension of the affine subspace that together define the feasible set. The number $\nnz$  is the total number of non-zero entries of the matrix $\tt A$ and cost vector $\tt c$ used to describe the problem in SeDuMi format.  These results show  problem size is often  significantly reduced and sparsity enhanced by our method.

\paragraph{2) DIMACS errors and distance to face}
We report a tuple $(e_1,\ldots,e_6)$ of DIMACS errors  \cite{mittelmann2003independent} for the original problem and reduced problem.
For instance, if the original problem has the form $(P-SDP)$, we solve it and report errors for $(P-SDP)$ and its dual $(D-SDP)$.  We then formulate  a reduced problem $(R/P-SDP)$ and report errors for $(R/P-SDP)$ and its dual $(R/D-SDP)$.  Finally, we report the distance $d_{face}$ (in norm induced by the trace inner-product) of the solution to the subspace spanned by the identified face. For instance, for an original SDP of the form $(P-SDP)$ with solution $y$, we report
\[
d_{face} = \| \Phi \left( \subspaceName(y)\right) - \subspaceName(y) \|_{F},
\]
where $\Phi : \mathbb{S}^n \rightarrow \mathbb{S}^n$ is the orthogonal projection map onto the mentioned subspace and $\| \cdot \|_F$ denotes the Frobenius norm. Note if the face equals $U \mathbb{S}^d_{+}U^T$ for $U$ with orthonormal columns, then 
\[
\Phi(X) = UU^TXUU^T.
\]
For original SDPs of the form $P-SDP$, the distance $d_{face}$  measures how well a solution $y$ satisfies the equations $V^T \subspaceName(y) U = 0$ and $V^T \subspaceName(y) V = 0$.  Note $d_{face}$ should be zero for exact solutions of both $(P-SDP)$ and $(R/P-SDP)$.

The reported errors show the reduced SDP can be solved just as accurately as the original in terms of DIMACS error.  They also show that by the measure $d_{face}$, solutions to the reduced SDP are significantly more accurate.  That $d_{face}$ is  larger for the original SDP  reflects the fact DIMACS error (a measure of backwards-error) can be  a poor measure of forwards-error when strict feasibility fails. (A phenomena observed in \cite{sturm2000error}.)

\paragraph{3) Reducing certificate error}
When one iteration of facial reduction is performed, we report
the minimum eigenvalue of the reducing certificate $S$ and a measure of the containment $\subspaceName \subseteq S^{\perp}$,
where $\subspaceName$ is the affine set of the SDP. (We restrict to one iteration to keep notation simple.) Specifically, for SDPs of the form $(P-SDP)$, we report a tuple $(| C \cdot S |, \max_i  | A_i \cdot S|,   \lambda_{\min}(S) )$
where the first two numbers measure containment of $\{ C-\sum^m_{i=1} y_i A_i :  y \in \mathbb{R}^m \}$ in $S^{\perp}$ and the last denotes
the minimum eigenvalue of $S$. For SDPs with feasible set defined by equations $A_i \cdot X = b_i$, we report $(| b^Ty |, \| S-\sum_i y_i A_i \|_F,   \lambda_{\min}(S) )$, where $(S,y)$   solves the appropriate variant of $(\star)$ described by equation \eqref{eq:perpPrim} and the discussion in Section~\ref{subsec:sdpfr}.   The reported errors show reducing certificates can be found with exceptional accuracy when
polyhedral approximations are used. 

\paragraph{4) Solve times}
For larger instances, we give solve times before and after reductions  and report  the total time  $t_{LPs}$  spent solving  LPs 
for reducing certificates. These solve times are reported for an Intel(R) Core(TM) i7-2600K CPU @ 3.40GHz machine with 16 gigabytes of RAM using the LP solver of MOSEK and the SDP solver SeDuMi called from MATLAB 2014a running Ubuntu.   For these instances,    solve time is  significantly reduced  and  the cost of solving LPs is negligible.

 \subsection{Lower bounds for optimal multi-period investment } \label{ex:boyd}
Our first example arises from SDP-based lower bounds of optimal multi-period investment strategies.  The strategies and specific SDP formulations are given in \cite{boyd2012performance}.
For each strategy, an SDP  produces a quadratic lower bound on the value function arising in the dynamic programming solution to the underlying optimization problem.
These bounds are produced using the $S$-procedure, an SDP-based method for showing emptiness
of sets defined by quadratic polynomials (see, e.g.,   \cite{boyd1994linear}). We report reductions  using diagonal  ($\D^d$)  approximations, DIMACs error, reducing certificate error, and solve time in  Tables~\ref{tab:invest1}-\ref{tab:invest3}. Scripts that generate the SDPs are found here (and require the package CVX \cite{grant2008cvx}):
\begin{center}
{\tt www.stanford.edu/\textasciitilde  boyd/papers/matlab/port\_opt\_bound/port\_opt\_code.tgz}
\end{center}

\newcommand{\vbarheavy}{1.5pt}

\begin{table}[tbp]
  \smaller
  \centering
   \subfloat[Original]{%
      \begin{tabular}{   | c | c| c| c|  }
      \hline
        Example &    $n$  &  $r$ &  $\nnz$  \\
      \hline
     {\tt long\_only}   &  $(91 \times 100, 30 \times 100)$ &  59095  &  853011   \\  \hline
     {\tt unconstrained}   &  $(121 \times 100, 30 \times 100)$ & 62095   &  874011   \\  \hline
     {\tt sector\_neutral}   &  $(121 \times 100, 30 \times 100)$ &  62392 &  1373000    \\  \hline
     {\tt leverage\_limit}   &  $(151 \times 100, 30 \times 100)$ & 68195 &  915993   \\  \hline
      \end{tabular}
    } 
    \subfloat[Reduced]{%
    \begin{tabular}{   | c | c| c| c|  }
      \hline
        Example &    $n$  &  $r$ &  $\nnz$  \\
      \hline
     {\tt long\_only}   &  $(61 \times 100, 30 \times 100)$ &  56095  &  832011   \\  \hline
     {\tt unconstrained}    &  $(61 \times 100, 30 \times 100)$ & 56095   &  840891   \\  \hline
     {\tt sector\_neutral} &  $(61 \times 100, 30 \times 100)$ &  56392 &  1342880    \\  \hline
     {\tt leverage\_limit}    &  $(61 \times 100, 30 \times 100)$ & 59195 &  873873   \\  \hline
      \end{tabular}
      \hspace{.5cm}%
    }  \caption{ Dimension $r$ of subspace and order $(n_1,\ldots,n_{200})$ of cone $\mathbb{S}_+^{n_{1}} \times \cdots \times \mathbb{S}_+^{n_{200}}$ describing feasible set. The column `$\nnz$' shows number of non-zero entries of SDP data matrices. } \label{tab:invest1}
  \subfloat[Original]{%
    \begin{tabular}{   | c | c | c| c|c|c|c!{\vrule width \vbarheavy} c| }
    \hline
     Example & $e_1$ &  $e_2$  & $e_3$  & $e_4$ &  $e_5$  & $e_6$  & $d_{face}$ \\
    \hline
{\tt long\_only} & 	4.3e-08 &	 0 &	 0 &	4.4e-11 &	-4.0e-05 &	-3.9e-05 &	2.3e-06  \\ 
 {\tt unconstrained} & 7.1e-08 &	 0 &	 0 &	1.1e-11 &	-2.5e-06 &	-1.6e-06 &	2.1e-05  \\ 
{\tt sector\_neutral} &	4.2e-07 &	 0 &	 0 &	1.1e-10 &	-1.4e-08 &	3.1e-05 &	1.6e-04  \\ 
{\tt leverage\_limit}	 & 7.3e-08 &	 0 &	 0 &	1.0e-11 &	-1.6e-06 &	-6.4e-07 &	1.2e-05  \\   \hline
    \end{tabular}
  } \\
  \subfloat[Reduced]{%
   \begin{tabular}{ | c | c | c| c|c|c|c!{\vrule width \vbarheavy} c|  }
    \hline
    Example & $e_1$ &  $e_2$  & $e_3$  & $e_4$ &  $e_5$  & $e_6$  & $d_{face}$ \\
    \hline
{\tt long\_only} & 	3.7e-08 &	 0 &	 0 &	1.5e-11 &	-6.8e-06 &	-5.8e-06 & 2.9e-17    \\ 
 {\tt unconstrained} & 	4.9e-08 &	 0 &	 0 &	1.2e-11 &	-4.9e-07 &	3.9e-07 &	3.1e-17    \\ 
{\tt sector\_neutral} &	3.5e-07 &	 0 &	 0 &	1.0e-10 &	-3.5e-08 &	3.0e-05 &	3.5e-17    \\ 
{\tt leverage\_limit}	 &	4.8e-08 &	 0 &	 0 &	9.7e-12 &	-1.1e-06 &	-2.9e-07 &	4.3e-14    \\  \hline
    \end{tabular}  
  } 
   \caption{DIMACS errors $e_i$ and distance $d_{face}$ to linear span of identified face.} \label{tab:invest2} 
   
  \vspace{.3cm}
     \begin{tabular}{ | c | c | c| c|c|c|c!{\vrule width \vbarheavy} c|  }
      \hline
      Example &  $| C \cdot S |$ & $\max_i  | A_i \cdot S|$  & $\lambda_{\min}(S)$    \\
      \hline
  {\tt long\_only} & 	0 &	 0 &	 0  	     \\ \hline
   {\tt unconstrained} & 	0 &	 0 &	 0     \\ \hline
  {\tt sector\_neutral} &	0 &	 0 &	 0       \\ \hline
  {\tt leverage\_limit}	 &	0 &	 0 &	 0  	   \\  \hline
      \end{tabular} 
      \caption{Reducing certificate error.  The first two columns measure containment of the SDP's affine subspace in the hyperplane $S^{\perp}$. The last denotes the minimum eigenvalue of the reducing certificate $S$.}

  \vspace{.3cm}
        \begin{tabular}{   | c | c| c|  c | c|c| }
        \hline
          Example &     Original     &  Reduced     & $\timeLP$ \\
        \hline
       {\tt long\_only}   &   651 &   613  & 0.33  \\  \hline
      {\tt unconstrained}   & 800  & 574  &  0.71 \\  \hline
      {\tt sector\_neutral}  &  760 & 496 & 0.70 \\  \hline
    {\tt leverage\_limit}  &  976 &  617  & 1.2   \\  \hline
        \end{tabular}
     \caption{Solve times (sec) for original and reduced SDPs.  The reduced SDP was formulated by solving
     LPs over diagonal approximations, i.e., by taking $\mathcal{C}(\setName) = \D^d$. These LPs took $\timeLP$ seconds to solve.
      } \label{tab:invest3} 
\end{table}

\subsection{Copositivity of quadratic forms} \label{ex:horn}

Our next example pertains to SDPs that demonstrate \emph{copositivity} of certain quadratic forms. A quadratic form
$x^T J x$ is copositive if and only if $x^T J x \ge 0$ for all $x$ in the non-negative orthant. 
Deciding copositivity is NP-hard, but a sufficient condition can be checked using sum-of-squares techniques and semidefinite programming,
as we now illustrate.

\paragraph{The Horn form} An example of a copositive polynomial is the Horn form $f(x) := x^T J x$, where
\[
 J = \left(\begin{array}{rrrrr} 1 & -1 & 1 & 1 & -1\\ -1 & 1 & -1 & 1 & 1\\ 1 & -1 & 1 & -1 & 1\\ 1 & 1 & -1 & 1 & -1\\ -1 & 1 & 1 & -1 & 1 \end{array}\right), \qquad  x = \left( \begin{array}{ccccc} x_1 & x_2 & x_3 & x_4 & x_5 \end{array} \right)^T.
\]
This polynomial, originally introduced by A. Horn, appeared previously in \cite{diananda1962non} \cite{quist1998copositive}. To see how copositivity can be demonstrated using SDP,  first 
note copositivity of $f(x)$ is equivalent to global non-negativity of   $f(z_1^2,z_2^2,z_3^2,z_4^2,z_5^2)$, where we have substituted each variable $x_i$
with the square of a new indeterminate $z^2_i$.  Next, note global non-negativity of the latter polynomial
 can be demonstrated by showing 
\begin{align} \label{eq:hornSOS}
 g(z) = \left(\sum^5_{i=1} z^2_i\right)  f(z_1^2,z_2^2,z_3^2,z_4^2,z_5^2)
\end{align}
is a sum-of-squares, which is equivalent to feasibility of a particular SDP over $\mathbb{S}^n_{+}$
where  $n={5+2 \choose 3}$, the number of    degree-three monomials in $5$ variables (see Chapter 3 of \cite{blekherman2013semidefinite} for details on constructing this SDP).

\paragraph{Generalized Horn forms}
The Horn form $f(x)$ generalizes to a family of copositive forms in  $n=3m+2$ variables ($m \ge 1)$:
\[
 B(x;m) = \left(\sum^{3m+2}_{i=1} x_i\right)^2 - 2 \sum^{3m+2}_{i=1} x_i \sum^m_{j=0} x_{i+3j+1},
\]
 where we let the subscript for the indeterminate $x$ wrap cyclically, i.e. $x_{r+n} = x_r$. This family was studied in \cite{baston1969extreme}, and the Horn form corresponds to the case $m=1$.  As with the Horn form, we can show copositivity of $B(x;m)$ by showing a polynomial analogous to \eqref{eq:hornSOS} is a sum-of-squares. We formulate SDPs that demonstrate copositivity of $B(x;m)$ in this way for each $m \in \{1,\ldots,5\}$. 
We report reductions  using diagonally-dominant ($\DD^d$)  approximations, DIMACs error, reducing certificate error, and solve time in  Tables~\ref{tab:horn1}-\ref{tab:horn3}. (Errors and solve time are omitted for $m>3$ since the SDPs are too large to solve.)

\begin{table}[tbp] \centering
   \subfloat[Original]{%
      \begin{tabular}{   | c | c| c| c|  }
      \hline
        Example &    $n$  &  $r$ &  $\nnz$  \\
      \hline
     $m=1$     & 35  &  420 &   1225  \\  \hline
     $m=2$     & 120  &  5544 &    14400  \\  \hline
     $m=3$     & 286   & 33033  &   81796     \\  \hline
     $m=4$     & 560   & 129948  &  313600      \\  \hline
     $m=5$     & 969   & 395352  &  938961      \\  \hline
      \end{tabular}
    } \qquad
    \subfloat[Reduced]{%
    \begin{tabular}{   | c | c| c| c|  }
      \hline
        Example &    $n$  &  $r$ &  $\nnz$  \\
      \hline
     $m=1$     & 25  & 165  &   1200  \\  \hline
     $m=2$     & 96  &  3132 &    14312  \\  \hline
     $m=3$     &  242  & 21879  &    81554    \\  \hline
     $m=4$     &  490  & 494143  &  313040      \\  \hline
     $m=5$     &  867  &  303399 &  937822      \\  \hline
      \end{tabular}
      \hspace{.5cm}%
    } \caption{ Dimension $r$ of subspace and order $n$ of cone $\mathbb{S}^n$ describing feasible set. The column `$\nnz$' shows number of non-zero entries of SDP data matrices. } \label{tab:horn1}
  \subfloat[Original]{%
   \begin{tabular}{ | c | c | c| c|c|c|c!{\vrule width \vbarheavy} c|  }
    \hline
    Example & $e_1$ &  $e_2$  & $e_3$  & $e_4$ &  $e_5$  & $e_6$  & $d_{face}$ \\
    \hline
     $m=1$  	& 8.63e-10  	& 0  	& 0  	& 9.99e-11  	& 1.99e-10  	& 3.14e-08  	& 8.47e-06  	 \\ \hline 
     $m=2$  	& 9.34e-09  	& 0  	& 0  	& 3.68e-10  	& 8.12e-10  	& 6.05e-07  	& 3.96e-05  	 \\ \hline 
     $m=3$  	& 1.87e-09  	& 0  	& 0  	& 1.01e-10  	& 1.97e-10  	& 4.16e-07  	& 3.83e-05  	 \\ \hline 
    \end{tabular} 
  } \\     
  \subfloat[Reduced]{%
    \begin{tabular}{   | c | c | c| c|c|c|c!{\vrule width \vbarheavy} c|  }
    \hline
     Example & $e_1$ &  $e_2$  & $e_3$  & $e_4$ &  $e_5$  & $e_6$  & $d_{face}$ \\
    \hline
     $m=1$  	& 7.82e-10  	& 0  	& 0  	& 2.42e-11  	& 7.64e-11  	& 2.04e-08  	& 9.27e-16  	 \\ \hline 
     $m=2$  	& 1.23e-09  	& 0  	& 0  	& 1.59e-10  	& 3.82e-10  	& 5.84e-08  	& 3.26e-16  	 \\ \hline 
     $m=3$  	& 4.00e-10  	& 0  	& 0  	& 7.08e-11  	& 2.25e-10  	& 7.93e-08  	& 7.48e-16  	 \\ \hline 
    \end{tabular} 
  } 
    \caption{DIMACs errors $e_i$ and distance $d_{face}$ to linear span of identified face.}  \label{tab:horn2}
 \vspace{.3cm}

          \begin{tabular}{ | c | c | c| c|c|c|c!{\vrule width \vbarheavy} c|  }
           \hline
           Example & $| b^Ty |$ &  $\| S-\sum_i y_i A_i \|_F$  & $\lambda_{\min}(S)$    \\
           \hline
      $m=1$  & 	0 &	 3.33e-16 &	 0  	     \\  \hline
       $m=2$ &   0 & 	 1.67e-16 &	 0 	    \\  \hline
       $m=3$  &	0 &	  -1.28e-15 &	 0   \\ \hline   
           \end{tabular} 
           \caption{Reducing certificate error.  The first two columns measure containment of the SDP's affine subspace in the hyperplane $S^{\perp}$. The last denotes the minimum eigenvalue of the reducing certificate $S$.}
       \vspace{.3cm}

        \begin{tabular}{   | c | c| c|  c | c|c| }
        \hline
          Example &     Original     &  Reduced     & $\timeLP$ \\
        \hline
     $m=1$   &    .81 &   .23   &  .047  \\  \hline
       $m=2$   &  11  &   9.2 &   .58 \\  \hline
      $m=3$   & 3900   & 3200   &  4.3  \\  \hline
        \end{tabular}
      \caption{Solve times (sec) for original and reduced SDPs.  The reduced SDP was formulated by solving
           LPs over diagonal approximations, i.e., by taking $\mathcal{C}(\setName) = \DD^d$. These LPs took $\timeLP$ seconds to solve.}  \label{tab:horn3}
\end{table}

\subsection{Lower bounds on completely positive rank} \label{ex:cprank}

A matrix $A \in \mathbb{S}^n$ is \emph{completely positive} (CP) if there exist $r$ non-negative vectors $v_i \in \mathbb{R}^n$ for which
\begin{align} \label{cp:pos}
 A = \sum^r_{i=1} v_i v_i^T.
\end{align}
The completely positive rank of $A$, denoted $\cprank A$, is the smallest  $r$ for which $A$ admits the decomposition \eqref{cp:pos}.
It follows trivially that
\begin{align*}
 \rank A \le \cprank A.
\end{align*}
In \cite{fawzi2014self}, Fawzi and the second author give an SDP formulation that improves this lower bound for a fixed matrix $A$.  This bound,  denoted $\tau^{sos}_{cp} (A)$ in \cite{fawzi2014self},
equals the optimal value of the following semidefinite program:
\begin{align*}
\begin{array}{crl}  \mbox{minimize } t \\ \mbox{subject to} \\ & \left( \begin{array}{cc} t & \vect A^T \\ \vect A & X \end{array} \right) \succeq 0 \\ 
                                                              & X_{ij,ij} \le A^2_{ij} & \forall i,j \in \{1,\ldots,n\}  \\
							      & X \preceq A \otimes A \\
							      & X_{ij,kl} = X_{il,jk} & \forall (1,1) \le (i,j) < (k,l) \le (n,n),
\end{array} 
\end{align*}
where $A \otimes A$ denotes the \emph{Kronecker product} and $\vect A$ denotes the $n^2 \times 1$ vector obtained by stacking the columns of $A$.
Here, the double subscript $ij$ indexes the $n^2$ rows (or columns) of $X$ 
and the inequalities on $(i,j)$   hold iff they hold element-wise (see  \cite{fawzi2014self} for  further clarification on this notation).

In this example, we formulate  SDPs as above for computing $\tau^{sos}_{cp}(Z)$,  $\tau^{sos}_{cp}(Z \otimes Z)$,   and   $\tau^{sos}_{cp}(Z\otimes Z \otimes Z)$,
where $Z$ is the completely positive matrix:
\[
Z = \left(\begin{array}{ccc} 4 & 0 & 1 \\ 0 & 4 & 1 \\ 1 & 1 & 3 \end{array}\right).
\]
Notice that since $Z$ is CP, the Kronecker products $Z \otimes Z$ and $Z \otimes Z \otimes Z$  are CP
(using the fact that $A \otimes B$ is $CP$ when $A$ and $B$ are CP \cite{bermanBook03}). Also notice that since
$Z$ contains zeros, the constraint $X_{ij,ij} \le Z^2_{ij}$ implies that $X$ 
has rows and columns identically zero; in other words, because $Z$ has elements equal to zero, the SDP for computing $\tau^{sos}_{cp}(Z)$
cannot have a strictly feasible solution.

To reduce the formulated SDPs, we first observe that each is actually a cone program over $\mathbb{R}^{n_1}_{+} \times \mathbb{S}^{n_2}_{+} \times \mathbb{S}^{n_3}_{+}$, i.e., each SDP
has a mix of linear inequalities and semidefinite constraints.   To find reductions,  we first treat the linear equalities as a semidefinite constraint on a diagonal matrix.
We report reductions  using diagonal ($\D^d$)  approximations, DIMACs error, reducing certificate error, and solve time in  Tables~\ref{tab:cprank1}-\ref{tab:cprank3}. (Solve times and  errors are omitted for $Z \otimes Z \otimes Z$ since the SDP is too large to solve.)

\begin{table}[tbp]
  \centering
   \subfloat[Complexity parameters - original]{%
      \begin{tabular}{   | c | c| c| c|  }
      \hline
        Example &    $n$  &  $r$ &  $\nnz$  \\
      \hline
     $Z$     &  $(9,10,9)$ & $37$ & $260$    \\  \hline
     $Z\otimes Z$     &  $(81,82,81)$ & $2026$ & $18344$     \\  \hline
     $Z \otimes Z \otimes Z$      &  $(729,730,729)$ & $142885$ & $1428692$      \\  \hline
      \end{tabular}
    }  \qquad
    \subfloat[Complexity parameters -  reduced]{%
    \begin{tabular}{   | c | c| c| c|  }
      \hline
        Example &    $n$  &  $r$ &  $\nnz$  \\
      \hline
     $Z$     &  $(7,8,9)$ & $20$ & $187$    \\  \hline
     $Z\otimes Z$     &  $(49,50,81)$ & $464$ & $8336$     \\  \hline
     $Z \otimes Z \otimes Z$      &  $(343,344,729)$ & $13262$ & $408403$      \\  \hline
      \end{tabular}
      \hspace{.5cm}%
    }   \caption{ Dimension $r$ of subspace and order $n$ of cone $\mathbb{R}^{n_1}_{+} \times \mathbb{S}^{n_1}_+ \times \mathbb{S}^{n_2}_+$ describing feasible set. The column `$\nnz$' shows number of non-zero entries of SDP data matrices. } \label{tab:cprank1}
  \subfloat[Original]{%
    \begin{tabular}{   | c | c | c| c|c|c|c!{\vrule width \vbarheavy} c|  }
    \hline
     Example & $e_1$ &  $e_2$  & $e_3$  & $e_4$ &  $e_5$  & $e_6$  & $d_{face}$ \\
    \hline
  $Z$  	   & 2.08e-11  	& 0  	   & 0  	& 1.09e-10 & -1.36e-09  	& -1.44e-09  	& 1.02e-05  	 \\ \hline 
  $Z \otimes Z$  	    & 6.58e-09 & 0  	& 0  	   & 1.72e-10  	& 2.82e-06  	& 2.68e-06  	& 2.42e-03 \\ \hline 
    \end{tabular}
  }     \\
  \subfloat[Reduced]{%
    \begin{tabular}{   | c | c | c| c|c|c|c!{\vrule width \vbarheavy} c|  }
    \hline
     Example & $e_1$ &  $e_2$  & $e_3$  & $e_4$ &  $e_5$  & $e_6$  & $d_{face}$ \\
    \hline
  $Z$  	& 1.27e-11  	& 0  	& 0  	& 7.87e-11  	& -1.54e-09  	& -1.59e-09  	& 0  	 \\ \hline 
  $Z \otimes Z$  	& 3.91e-08  	& 0  	& 0  	& 0  	& 8.50e-06  	& 7.56e-06  	& 0  	 \\ \hline 
    \end{tabular}
  } \caption{DIMACs errors $e_i$ and distance $d_{face}$ to linear span of identified face.}\label{tab:cprank2} 
\vspace{.35cm}

       \begin{tabular}{ | c | c | c| c|c|c|c!{\vrule width \vbarheavy} c|  }
        \hline
        Example &    $| C \cdot S |$ & $\max_i  | A_i \cdot S|$  & $\lambda_{\min}(S)$    \\
        \hline
   $Z$  & 	0 &	 0 &	 0  	     \\ \hline 
    $Z \otimes Z$   & 	0 &	 0 &	 0     \\ \hline 
        \end{tabular} 
  \caption{Reducing certificate error.  The first two columns measure containment of the SDP's affine subspace in the hyperplane $S^{\perp}$. The last denotes the minimum eigenvalue of the reducing certificate $S$.}
\vspace{.35cm}
        \begin{tabular}{   | c | c| c|  c | c|c| }
        \hline
          Example &     Original     &  Reduced     & $\timeLP$ \\
        \hline
     $Z$   &   .4 &   .7  & .0084  \\  \hline
      $Z \otimes Z$   & 131  & 10.5  &  .016 \\  \hline
        \end{tabular}
     \caption{Solve times (sec) for original and reduced SDPs.  The reduced SDP was formulated by solving
          LPs over diagonal approximations, i.e., by taking $\mathcal{C}(\setName) = \D^d$. These LPs took $\timeLP$ seconds to solve. } \label{tab:cprank3}
\end{table}

\subsection{Lyapunov Analysis of a Hybrid Dynamical System}\label{ex:lyap}
The next example arises from  SDP-based stability analysis of a \emph{rimless wheel}, a hybrid dynamical system and simple model for walking robots  studied in
\cite{posa13a} by  Posa, Tobenkin, and Tedrake. The SDP  includes several coupled semidefinite constraints that impose Lyapunov-like stability conditions accounting for Coulomb friction and the contact dynamics of the rimless wheel.    We report reductions using diagonally-dominant $(\DD^d)$ and diagonal ($\D^d$)  approximations, DIMACs error, and solve time in  Tables~\ref{tab:lyap1}-\ref{tab:lyap3}. (Reducing certificate error is omitted since multiple facial reduction iterations were performed.)

\begin{table} 
\centering
   \begin{tabular}{   | c | c| c|c| }
    \hline
    Problem & $n$ &  $r$ & $\nnz$ \\
    \hline
  Original  	& $(6  	, 108  	, 11  	, 11  	, 11  	, 11  	, 11  	, 11  	, 11  	, 11  	, 11  	, 11)$  	& 4334  	& 16864  	 \\ \hline 
    Reduced, $\mathcal{C}(\setName) =\D^d$   	& $(6  	, 56  	, 11  	, 1  	, 1  	, 0  	, 11  	, 1  	, 1  	, 0  	, 11  	, 11)$  	& 1138  	& 6661  	 \\ \hline 
   Reduced, $\mathcal{C}(\setName) = \DD^d$    	& $(6  	, 34  	, 8  	, 1  	, 1  	, 0  	, 8  	, 1  	, 1  	, 0  	, 9  	, 7)$  	& 452  	& 4007  	 \\ \hline 
    \end{tabular} \caption{The feasible set is an $r$-dimensional subspace intersected with the cone $\mathbb{S}^{n_1}_{+}\times \mathbb{S}^{n_2}_{+} \cdots  \times \mathbb{S}^{n_{12}}_{+}$.} \label{tab:lyap1}
\vspace{.4cm}
   \begin{tabular}{ | c|c|c|c|c|c|c!{\vrule width \vbarheavy} c| }
    \hline
    Problem & $e_1$ &  $e_2$  & $e_3$ & $e_4$ & $e_5$ & $e_6$ & $d_{face}$ \\
    \hline
  Original  	& 2.56e-07  	& 0  	& 0  	& 2.48e-10  	& 9.78e-08  	& 1.70e-05  	& 8.99e-02  	 \\ \hline 
          Reduced, $\mathcal{C}(\setName) =\D^d$    	& 2.77e-08  	& 0  	& 0  	& 0  	& 1.76e-08  	& 8.29e-06  	& 0  	 \\ \hline 
      Reduced, $\mathcal{C}(\setName) = \DD^d$ 	& 6.65e-08  	& 0  	& 0  	& 0  	& 4.29e-08  	& 1.20e-05  	& 3.82e-15  	 \\ \hline 
    \end{tabular}  \caption{DIMACS error bounds $e_i$ and distance $d_{face}$ to the linear span of  identified face.} \label{tab:lyap2}
\vspace{.4cm}  
   \begin{tabular}{   | c | c| c| }
    \hline
    Problem & Solve time & $t_{LPs}$ \\
    \hline
  Original  	&  111  	& --	 \\ \hline 
    Reduced, $\mathcal{C}(\setName) =\D^d$   	& 5	&   .05	  	 \\ \hline 
   Reduced, $\mathcal{C}(\setName) = \DD^d$    	& 1.8 	&   0.82   	 \\ \hline 
    \end{tabular} 
  \caption{
Solve times (sec) for original and reduced SDPs.  The reduced SDP was formulated by solving  LPs over the indicated approximation ($\mathcal{C}(\setName) = \D^d$ or $\mathcal{C}(\setName) = \DD^d$) which took $\timeLP$ seconds to solve.} \label{tab:lyap3}
\end{table}

\subsection{Multi-affine polynomials, matroids, and the half-plane property} \label{ex:matroids}

A multivariate polynomial $f(z) : \mathbb{C}^n \rightarrow \mathbb{C}$ has the \emph{half-plane} property if it is non-zero when each variable $z_i$ has positive real part.
A polynomial is \emph{multi-affine} if  each indeterminate is raised to at most the first power.
As proven in \cite{choe2004homogeneous}, if a multi-affine, homogeneous polynomial with unit coefficients has the half-plane property, it is the \emph{basis generating polynomial} of a \emph{matroid}. In this section, we reduce SDPs 
that arise in the study of the converse question:  given a matroid, does its basis generating polynomial have the half-plane property?  Or more precisely, given a rank-$r$ matroid $M$ (over the ground-set $\{1,\ldots,n\}$) with  set of bases $B(M)$,   does 
the multi-affine, degree-$r$ polynomial  
\begin{align} \label{eq:matroid}
  f_M(z_1,\ldots,z_n) := \sum_{ \substack{ \{i_1,i_2,\ldots,i_r \} \\ \in B(M) } } z_{i_1} z_{i_2} \cdots z_{i_r}
\end{align}
 have the half-plane property?  

\paragraph{The role of polynomial non-negativity}
 This converse question is related to global non-negativity of so-called \emph{Rayleigh differences} of $f_M(z)$, which are polynomials over $\mathbb{R}^n$  defined for each $\{i,j\} \subset \{ 1, \ldots, n \}$
as follows:
\[
 \Delta_{ij} f_M(x) := \frac{\partial f_M}{\partial z_i}(x) \frac{\partial f_M}{\partial z_j}(x)-\frac{\partial^2 f_M}{\partial z_i \partial z_j}(x) \cdot f_M(x).
\]
 A theorem of Br{\"a}nd{\'e}n  \cite{branden2007polynomials}  states  $f_M(z)$ has the half-plane property if and only if all of ${n \choose 2}$ Rayleigh differences are globally non-negative, i.e., $\Delta_{ij} f_M(x) \ge 0$ for all $x \in \mathbb{R}^n$.  
 An equivalent criterion, stated in terms of global non-negativity of   a single Rayleigh difference (and so-called \emph{contractions} and \emph{deletions} of $M$), appears in \cite{wagner2009criterion}.

\paragraph{The role of semidefinite programming}
Since semidefinite programming can  demonstrate a given polynomial is a sum-of-squares, it is a natural tool for proving
a given Rayleigh difference $\Delta_{ij} f_M(x)$ is globally non-negative. In this section,
we  formulate and then apply our reduction technique to SDPs  that test the sum-of-squares condition for various $\Delta_{ij} f_M(x)$  and various matroids $M$. 
As is standard, the SDPs are formulated using the  set of monomial exponents in $\frac{1}{2} \mathcal{N} (\Delta_{ij} f_M) \cap \mathbb{N}^{n}$,
 where $\mathcal{N} (\Delta_{ij}f_M)$ denotes the \emph{Newton polytope} of $\Delta_{ij}f_M$ (see Chapter 3 of \cite{blekherman2013semidefinite} for details on this formulation).

We report reductions  using diagonally-dominant $\DD^d$ approximations, DIMACs error, and reducing certificate error in  Tables~\ref{tab:Matroids1}-\ref{tab:Matroids3}. (Solve time is omitted since the original SDPs are small.)  We now elaborate on each matroid in these tables.

\subsubsection{Various matroids  with the half-plane property} \label{sec:weimat}

The first set of matroids were studied by \citet{wagner2009criterion}.  Specifically, \citet{wagner2009criterion} demonstrate that  $\Delta_{ij} f_{M}$  (for specific $\{i,j\}$) is a sum-of-squares
for  matroids $M$  they denote $\mathcal{F}^{-4}_7$, $\mathcal{W}^{3+}$, $\mathcal{W}^{3}+e$,
$\mathcal{P}^{'}_7$, $n\mathcal{P} \setminus 1$, $n\mathcal{P} \setminus 9$, and $\mathcal{V}_8$.
(We refer the reader to  \cite{wagner2009criterion} for definitions of these matroids and the explicit polynomials $\Delta_{ij} f_{M}$.)  Note Wagner and Wei demonstrate each sum-of-squares condition via ad-hoc construction, instead of by solving an SDP.

Notice from Table~\ref{tab:Matroids1} that for matroids $\mathcal{W}^{3+}$, $\mathcal{W}^{3}+e$, $\mathcal{P}^{'}_7$, $n\mathcal{P} \setminus 1$ and $n\mathcal{P} \setminus 9$, the reduced SDP is described by a zero-dimensional affine subspace. In other words, the SDP demonstrating the sum-of-squares condition has a feasible set containing a single point.

\begin{table}[tbp]
  \small
  \centering
   \subfloat[Original]{%
    \begin{tabular}{   | c |c| c| c| c|  }
      \hline
        Matroid & $\{i,j\}$ &  $n$  &  $r$ &  $\nnz$  \\
      \hline
$\mathcal{F}^{-4}_7$ & $\{1,2\}$  	& 8  	& 5  	& 64  	 \\ \hline 
$\mathcal{W}^{3+}$ & $\{1,2\}$  	& 8  	& 5  	& 64  	 \\ \hline 
$\mathcal{W}^{3}+e$ & $\{1,2\}$  	& 9  	& 7  	& 81  	 \\ \hline 
$\mathcal{P}^{'}_7$ &  $\{1,2\}$  	& 8  	& 4  	& 64  	 \\ \hline 
$n\mathcal{P}  \setminus 1$ & $\{2,4\}$  	& 12  	& 14  	& 144  	 \\ \hline 
$n\mathcal{P}  \setminus 9$ & $\{1,2\}$  	& 12  	& 14  	& 144  	 \\ \hline 
   $\mathcal{V}_8$ & $\{1,2\}$  	& 16  	& 33  	& 256  	 \\ \hline 
$\mathcal{V}_{10}$ & $\{3,4\}$  	& 52  	& 657  	& 2704  	 \\ \hline 
      \end{tabular}
    } \hspace{1cm}
    \subfloat[Reduced]{%
    \begin{tabular}{   | c |c| c| c| c|  }
      \hline
        Matroid & $\{i,j\}$ &  $n$  &  $r$ &  $\nnz$  \\
      \hline
$\mathcal{F}^{-4}_7$ & $\{1,2\}$  	& 5  	& 1  	& 25  	 \\ \hline 
$\mathcal{W}^{3+}$ & $\{1,2\}$  	& 3  	& 0  	& 9  	 \\ \hline 
$\mathcal{W}^{3}+e$ & $\{1,2\}$  	& 5  	& 0  	& 27  	 \\ \hline 
$\mathcal{P}^{'}_7$ &  $\{1,2\}$  	& 4  	& 0  	& 16  	 \\ \hline 
$n\mathcal{P}  \setminus 1$ & $\{2,4\}$  	& 6  	& 0  	& 40  	 \\ \hline 
$n\mathcal{P}  \setminus 9$ & $\{1,2\}$  	& 5  	& 0  	& 27  	 \\ \hline 
   $\mathcal{V}_8$ & $\{1,2\}$  	& 13  	& 17  	& 185  	 \\ \hline 
$\mathcal{V}_{10}$ & $\{3,4\}$  	& 41  	& 327  	& 2087  	 \\ \hline  
      \end{tabular}
      \hspace{.5cm}%
    } \caption{ Dimension $r$ of subspace and order $n$ of cone $\mathbb{S}_+^n$ describing feasible set. The column `$\nnz$' shows number of non-zero entries of SDP data matrices. } \label{tab:Matroids1}
  \subfloat[Original]{%
    \begin{tabular}{   |c| c | c | c| c|c|c|c!{\vrule width \vbarheavy} c| }
    \hline
     Matroid & $\{i,j\}$  & $e_1$ &  $e_2$  & $e_3$  & $e_4$ &  $e_5$  & $e_6$  & $d_{face}$ \\
    \hline
$\mathcal{F}^{-4}_7$ & $\{1,2\}$  	& 2.10e-12  	& 0  	& 0  	& 9.32e-13  	& 7.50e-12  	& 2.06e-11  	& 3.25e-12  	 \\ \hline 
$\mathcal{W}^{3+}$ & $\{1,2\}$  	& 7.57e-11  	& 0  	& 0  	& 3.52e-11  	& 6.47e-11  	& 7.29e-10  	& 4.41e-06  	 \\ \hline 
$\mathcal{W}^{3}+e$ & $\{1,2\}$  	& 8.14e-09  	& 0  	& 0  	& 1.05e-11  	& 2.21e-11  	& 4.04e-10  	& 1.74e-06  	 \\ \hline 
$\mathcal{P}^{'}_7$ &  $\{1,2\}$  	& 2.31e-10  	& 0  	& 0  	& 8.55e-10  	& 1.39e-09  	& 2.19e-09  	& 6.54e-06  	 \\ \hline 
$n\mathcal{P}  \setminus 1$ & $\{2,4\}$  	& 9.04e-10  	& 0  	& 0  	& 1.29e-10  	& 2.40e-10  	& 9.31e-09  	& 4.51e-06  	 \\ \hline 
$n\mathcal{P}  \setminus 9$ & $\{1,2\}$  	& 2.45e-09  	& 0  	& 0  	& 3.54e-10  	& 7.07e-10  	& 1.87e-08  	& 4.15e-08  	 \\ \hline 
   $\mathcal{V}_8$ & $\{1,2\}$  	& 5.29e-11  	& 0  	& 0  	& 8.32e-11  	& 1.78e-10  	& 6.49e-10  	& 4.19e-06  	 \\ \hline 
$\mathcal{V}_{10}$ & $\{3,4\}$  	& 3.64e-11  	& 0  	& 0  	& 1.40e-09  	& 2.73e-09  	& 9.78e-09  	& 5.37e-06  	 \\ \hline 
    \end{tabular}
  } \\
  \subfloat[Reduced]{%
    \begin{tabular}{   |c| c | c | c| c|c|c|c!{\vrule width \vbarheavy} c|}
    \hline
     Matroid & $\{i,j\}$ & $e_1$ &  $e_2$  & $e_3$  & $e_4$ &  $e_5$  & $e_6$  & $d_{face}$ \\
    \hline
$\mathcal{F}^{-4}_7$ & $\{1,2\}$  	& 4.65e-10  	& 0  	& 0  	& 2.16e-08  	& 6.61e-08  	& 6.86e-08  	& 0  	 \\ \hline 
$\mathcal{W}^{3+}$ & $\{1,2\}$  	& 1.35e-15  	& 0  	& 0  	& 2.91e-12  	& 5.55e-12  	& 5.55e-12  	& 4.57e-16  	 \\ \hline 
$\mathcal{W}^{3}+e$ & $\{1,2\}$  	& 8.86e-15  	& 0  	& 0  	& 1.27e-11  	& 2.03e-11  	& 2.03e-11  	& 4.62e-16  	 \\ \hline 
$\mathcal{P}^{'}_7$ &  $\{1,2\}$  	& 6.91e-16  	& 0  	& 0  	& 1.13e-11  	& 1.62e-11  	& 1.62e-11  	& 4.73e-16  	 \\ \hline 
$n\mathcal{P}  \setminus 1$ & $\{2,4\}$  	& 5.27e-11  	& 0  	& 0  	& 1.89e-08  	& 3.67e-08  	& 3.75e-08  	& 3.46e-16  	 \\ \hline 
$n\mathcal{P}  \setminus 9$ & $\{1,2\}$  	& 1.18e-10  	& 0  	& 0  	& 3.18e-08  	& 4.29e-08  	& 4.35e-08  	& 3.44e-16  	 \\ \hline 
   $\mathcal{V}_8$ & $\{1,2\}$  	& 1.43e-11  	& 0  	& 0  	& 4.35e-11  	& 7.77e-11  	& 1.97e-10  	& 0  	 \\ \hline 
$\mathcal{V}_{10}$ & $\{3,4\}$  	& 2.90e-11  	& 0  	& 0  	& 1.11e-09  	& 2.09e-09  	& 3.15e-09  	& 5.15e-17  	 \\ \hline 
    \end{tabular}
   }  

    \caption{DIMACs errors $e_i$ and distance $d_{face}$ to linear span of identified face.}\label{tab:Matroids2} \vspace{.3cm}
    
    \begin{tabular}{   |c| c | c | c| c|c|c|c!{\vrule width \vbarheavy} c|}
           \hline
            Matroid & $\{i,j\}$ & $| b^Ty |$ &  $\| S-\sum_i y_i A_i \|_F$  & $\lambda_{\min}(S)$\\
           \hline
    $\mathcal{F}^{-4}_7$ & $\{1,2\}$  	& 0  	& 0  	& 0  	 \\ \hline 
    $\mathcal{W}^{3+}$ & $\{1,2\}$  	& 2.22e-16  	& 0  	& 0  	 \\ \hline 
    $\mathcal{W}^{3}+e$ & $\{1,2\}$  	& 0  	& 0  	& 0  	 \\ \hline 
    $\mathcal{P}^{'}_7$ &  $\{1,2\}$  	& 0  	& 0  	& 0  	 \\ \hline 
    $n\mathcal{P}  \setminus 1$ & $\{2,4\}$  	& 6.66e-16  	& 0  	& 0  	 \\ \hline 
    $n\mathcal{P}  \setminus 9$ & $\{1,2\}$  	& 0  	& 0  	& 0  	 \\ \hline 
       $\mathcal{V}_8$ & $\{1,2\}$  	& 0  	& 0  	& 0  	 \\ \hline 
    $\mathcal{V}_{10}$ & $\{3,4\}$  	& 1.78e-15  	& 0  	& 0  	 \\ \hline 
           \end{tabular} 
              \caption{Reducing certificate error.  The first two columns measure containment of the SDP's affine subspace in the hyperplane $S^{\perp}$. The last denotes the minimum eigenvalue of the reducing certificate $S$. }\label{tab:Matroids3}

\end{table}

\subsubsection{Extended V{\'a}mos matroid} \label{sec:extvamos}
The other matroid considered was studied by Burton, Vinzant, and Youm in  \cite{Vinzant14}. There, the authors use semidefinite programming to show
$\Delta_{ij} f_{\mathcal{V}_{10}}$ is a sum-of-squares for a specific $\{i,j\}$, where
$\mathcal{V}_{10}$ denotes the  \emph{extended V{\'a}mos matroid}  defined over the ground set $\{1,\ldots,10\}$.
The bases of $\mathcal{V}_{10}$  are all cardinality-four subsets of $\{1,\ldots,10\}$   excluding 
\[
 \{1,2,6,7 \}, \{1,3,6,8 \}, \{1,4,6,9 \}, \{1,5,6,10 \}, \{2,3,7,8 \},  \{3,4,8,9 \},  \mbox{ and }\{4,5,9,10 \}.
\]
From these bases, we construct $f_{\mathcal{V}_{10}}$ via \eqref{eq:matroid} and formulate an SDP demonstrating $\Delta_{34} f_{\mathcal{V}_{10}}$  
is a sum-of-squares (as was done in \cite{Vinzant14}).

\subsection{Facial Reduction Benchmark Problems}
In \cite{cheung2011preprocessing}, Cheung, Schurr, and Wolkowicz developed a facial reduction procedure for identifying faces in a numerically stable manner.
They also created a set of benchmark problems for testing their method. These problem instances are available at the URL below: 
\begin{center}
{\tt http://www.math.uwaterloo.ca/\textasciitilde hwolkowi/henry/reports/SDPinstances.tar}.
\end{center}
Each problem is a primal-dual pair hand-crafted so that both the primal and dual have no strictly feasible solution.  
We apply our technique to each primal problem and each dual problem individually, using diagonally-dominant ($\DD^{d}$) approximations.  Results are shown in Table~\ref{tab:benchmark}. Since some of the examples have duality gaps, we do not show DIMACs errors  nor we do show solve time given the small sizes. We also omit reducing certificate error since multiple facial reduction iterations were performed.

\begin{table}
\begin{center}
  \begin{tabular}{ |c | c | c|c|c|c|c|  }
    \hline
    Example & \shortstack{ \\ Original \\ Primal \\ $n;r$} &  \shortstack{ \\ Reduced \\ Primal \\ $n;r$} &   \shortstack{ \\ Original 
\\ Dual \\ $n;r$} &  \shortstack{ \\ Reduced \\ Dual \\ $n;r$}   \\
    \hline
   {\tt Example1} & $3;4$ & $2;2$ &   $3;2$ & $1;1$  \\  \hline
    {\tt Example2} & $3;4$ & $2;2$ &   $3;2$ & $2;1$   \\  \hline
    {\tt Example3} & $3;2$ & $2;2$ &  $3;4$ & $2;2$   \\  \hline
    {\tt Example4} & $3;3$ & $1;0$ &   $3;3$ & $1;1$   \\  \hline
    {\tt Example5} & $10;50$ & $10;50$ &   $10;5$ & $10;5$    \\  \hline
    {\tt Example6} & $8;28$ & $5;11$ &   $8;8$ & $4;4$   \\  \hline
    {\tt Example7} & $5;12$ & $4;8$ &  $5;3$ & $1;1$  \\  \hline
    {\tt Example9a}  & $100;4950$ &  $1;0$ &   $100;100$ & $1;1$    \\  \hline
    {\tt Example9b}  & $20;190$ & $1;0$ &   $20;20$ & $1;1$  \\
\hline
  \end{tabular}
\end{center}
\caption{Complexity parameters for the primal-dual SDP pairs given in \cite{cheung2011preprocessing}. The feasible set of each SDP is an $r$-dimensional
subspace intersected with the cone $\mathbb{S}^n_{+}$.  To formulate each reduced SDP, a face was identified by solving  LPs  over diagonally-dominant approximations ($\DD^d$). These LPs took  (in total) $\timeLP$ seconds to solve.   } \label{tab:benchmark}
\end{table}

\subsection{Difficult SDPs arising in polynomial non-negativity} \label{sec:exStrange}
In \cite{waki2012generate} and \cite{waki2012strange}, Waki et al. study two sets of SDPs
that are difficult to solve. For one set of SDPs,  SeDuMi fails to find certificates of infeasibility \cite{waki2012generate}.
 For the other set,   SeDuMi reports an incorrect optimal value \cite{waki2012strange}. The sets of SDPs are available at:
\begin{center}
{ \tt https://sites.google.com/site/hayatowaki/Home/difficult-sdp-problems}.
\end{center}
It turns out for each primal-dual pair in these sets, the problem defined by equations $A_i\cdot X = b_i$    is not strictly feasible. We apply our technique 
to both sets of SDPs using diagonal approximations $\D^d$ 
and arrive at SDPs that are more easily solved.  In particular, certificates of infeasibility are found for the SDPs in \cite{waki2012generate} and correct 
optimal values  are found for the SDPs in \cite{waki2012strange} by solving the reduced SDPs with SeDuMi. Problem size reductions are shown in Table~\ref{tab:strange1} and Table~\ref{tab:strange2}.  We omit solve time comparisons and DIMACs errors since the reduced problem is a trivial SDP in each case. We omit reducing certificate error since multiple facial reduction iterations were performed.

 \begin{table} 
\begin{center}
    \begin{tabular}{   | c | c| c|  c |c| }
    \hline
      Example &   \shortstack{\\$n;r$ \\ Original} & \shortstack{\\ $n;r$  \\ Reduced }      \\
    \hline
   {\tt CompactDim2R1}   & 3;4 &  1;1       \\  \hline
   {\tt CompactDim2R2}   &(6,3,3,3); 25  & (1,0,1,1); 1    \\       \hline
   {\tt CompactDim2R3}    & (10,6,6,6); 91 &  (1,0,1,1); 1    \\       \hline
   {\tt CompactDim2R4}    & (15,10,10,10); 241 & (1,0,1,1); 1  \\       \hline
   {\tt CompactDim2R5}    & (21,15,15,15); 526 &  (1,0,1,1); 1   \\       \hline
   {\tt CompactDim2R6}    & (28,21,21,21); 1009  &  (1,0,1,1); 1 \\       \hline
   {\tt CompactDim2R7}    &  (36,28,28,28); 1765 & (1,0,1,1); 1    \\       \hline
   {\tt CompactDim2R8}    & (45,36,36,36); 2881  & (1,0,1,1); 1   \\       \hline
   {\tt CompactDim2R9}     & (55,45,45,45); 4456  & (1,0,1,1); 1  \\       \hline
   {\tt CompactDim2R10}     & (66,55,55,55); 6601  &  (1,0,1,1); 1   \\       \hline
    \end{tabular} \caption{
Complexity parameters for weakly-infeasible SDPs studied in \cite{waki2012generate}. The feasible set of each SDP is an $r$-dimensional
subspace intersected with the cone $\mathbb{S}^{n}_{+}$.   To formulate each reduced SDP, a face was identified by solving  LPs defined by diagonal approximations ($\D^d$).  These LPs  took (in total) $\timeLP$ seconds to solve.  
} \label{tab:strange1}
\end{center}

\end{table}

 \begin{table} 
\begin{center}
    \begin{tabular}{   | c | c| c| c| c | }
    \hline
    Example &  \shortstack{\\$n;r$ \\ Original} & \shortstack{\\$n;r$ \\ Reduced}    & \shortstack{ \\Optimal Value \\ Reduced }  \\
    \hline
  {\tt unboundDim1R2}     & (3,2,2); 8  & (1,1,0); 1    & 1.080478e-13 \\       \hline
  {\tt unboundDim1R3}     & (4,3,3); 16  & (1,1,0); 1    &1.080478e-13 \\       \hline
   {\tt unboundDim1R4}     & (5,4,4); 27  & (1,1,0); 1    & 1.080478e-13 \\       \hline
   {\tt unboundDim1R5}     & (6,5,5); 41 &  (1,1,0); 1  & 1.080478e-13 \\       \hline
   {\tt unboundDim1R6}     & (7,6,6); 58  &  (1,1,0); 1    & 1.080478e-13 \\       \hline
   {\tt unboundDim1R7}     & (8,7,7); 78  &  (1,1,0); 1    & 1.080478e-13 \\       \hline
   {\tt unboundDim1R8}     & (9,8,8); 101  &  (1,1,0); 1   & 1.080478e-13 \\       \hline
   {\tt unboundDim1R9}     & (10,9,9); 127  &  (1,1,0); 1   & 1.080478e-13 \\       \hline
   {\tt unboundDim1R10}     & (11,11,10); 156  &  (1,1,0); 1    & 1.080478e-13 \\       \hline
    \end{tabular} \caption{
Complexity parameters for the SDPs in \cite{waki2012strange}. The feasible set of each SDP is an $r$-dimensional
subspace intersected with the cone $\mathbb{S}^{n}_{+}$.   To formulate each reduced SDP, a face was identified by solving  LPs defined by diagonal approximations ($\D^d$).  These LPs  took (in total) $\timeLP$ seconds to solve. 
For these examples, SeDuMi incorrectly returns an optimal value of one for the original problem. The optimal value returned for the reduced problem is very near the correct optimal value of zero. 
} \label{tab:strange2}
\end{center}
\end{table}

\subsection{DIMACS Controller Design Problems} \label{ex:dimacs}

Our final examples are the controller design problems {\tt hinf12} and  {\tt hinf13} of the DIMACS library \cite{pataki1999dimacs}---which evidently are SDPs in the library
with no strictly feasible solution. Results are shown in Tables~\ref{tab:dimacs1}-\ref{tab:dimacs3}, where we apply facial reduction to the primal problem of both SDPs (using $\DD^d$ for {\tt hinf12} and $\SDD^d$ for {\tt hinf13}).  As observed in \cite{mittelmann2003independent}, these problem instances are extremely difficult for  SDP solvers. For purposes of comparison, we therefore report DIMACS errors for both SeDuMi and SDPT3 \cite{toh1999sdpt3}. Solution times are omitted given the small sizes of these SDPs.

\begin{table}[tbp]
  \small
  \centering
   \subfloat[Original]{%
    \begin{tabular}{   | c |c| c| c| c|  }
      \hline
        Matroid  &  $n$  &  $r$ &  $\nnz$  \\
      \hline
       {\tt hinf12}  	& (6,  6, 12)  	& 77  	& 990  	 \\ \hline 
       {\tt hinf13}  	& (7, 9,   14)  	& 121  	& 2559  	 \\ \hline 
      \end{tabular}
    } \hspace{1cm}
    \subfloat[Reduced]{%
    \begin{tabular}{   | c |c| c| c| c|  }
      \hline
       Problem  &  $n$  &  $r$ &  $\nnz$  \\
      \hline
      {\tt hinf12}    & (6,  2,  6)  	& 23  	& 583  	 \\ \hline 
     {\tt hinf13}  	  &  (1, 9,   7)  	& 45   	& 1465  	 \\ \hline 
      \end{tabular}
      \hspace{.5cm}%
    } \caption{ Dimension $r$ of subspace and order $n$ of cone $\mathbb{S}_+^n$ describing feasible set. The column `$\nnz$' shows number of non-zero entries of SDP data matrices. For ${\tt hinf12}$, we used $\DD^d$. For  ${\tt hinf13}$, we used $\SDD^d$. } \label{tab:dimacs1}
  \subfloat[Original]{%
    \begin{tabular}{   |c| c | c | c| c|c|c|c!{\vrule width \vbarheavy} c| }
    \hline
     Problem & $e_1$ &  $e_2$  & $e_3$  & $e_4$ &  $e_5$  & $e_6$  & $d_{face}$ \\ \hline
 {\tt hinf12/sedumi}  	& 5.04e-09  	& 0  	& 0  	& 0  	& -1.55e-02  	& 2.23e-01  	& 1.17e-08  	 \\ \hline 
 {\tt hinf13/sedumi}  	& 6.21e-05  	& 0  	& 0  	& 2.63e-06  	& -3.68e-03  	& 2.30e-02  	& 1.00e+00$^\dagger$	 \\ \hline
 {\tt hinf12/sdpt3}  	& 1.67e-11  	& 0   	& 1.72e-05  	& 0   	& -1.72e-06  	& 2.36e-05  	& 3.81e-12  	 \\ \hline 
 {\tt hinf13/sdpt3}  	& 9.97e-06  	& 0   	& 5.73e-07  	& 0  	& -2.35e-04  	& 1.94e-04  	& 1.43e-02  	 \\ \hline 
    \end{tabular}
  } \\
  \subfloat[Reduced]{%
    \begin{tabular}{   |c| c | c | c| c|c|c|c!{\vrule width \vbarheavy} c|}
    \hline
     Problem & $e_1$ &  $e_2$  & $e_3$  & $e_4$ &  $e_5$  & $e_6$  & $d_{face}$ \\
    \hline
       {\tt hinf12/sedumi}  	& 4.99e-09  	& 0  	& 0  	& 0  	& -5.62e-02  	& 2.82e-01  	& 0  	 \\ \hline 
       {\tt hinf13/sedumi}  	& 6.39e-05  	& 0  	& 0  	& 1.51e-06  	& -2.76e-04  	& 1.93e-03  	& 0  	 \\ \hline
       {\tt hinf12/sdpt3}  	& 1.58e-11  	& 0  	& 3.18e-06  	& 0  	& -2.06e-06  	& 3.33e-05  	& 0   	 \\ \hline 
       {\tt hinf13/sdpt3}  	& 3.84e-05  	& 0  	& 7.09e-08  	& 0  	& -6.61e-04  	& 1.07e-05  	& 0   	 \\ \hline 
    \end{tabular} 
   }  \caption{DIMACs errors $e_i$ and distance $d_{face}$ to linear span of identified face. Normalized by solution norm, the outlier,  marked $^\dagger$, equals $2.53e-04$.   }\label{tab:dimacs2}  \vspace{.3cm}  
             \begin{tabular}{ | c | c | c| c|c|c|c!{\vrule width \vbarheavy} c|  }
              \hline
              Example & $| b^Ty |$ &  $\| S-\sum_i y_i A_i \|_F$  & $\lambda_{\min}(S)$    \\
              \hline
          {\tt hinf12}  & 	0 &	 0 &	 0  	     \\  \hline
         {\tt hinf13} &   0 & 	  8.31e-10 &	 0 	    \\  \hline
              \end{tabular} 
              \caption{Reducing certificate error.  The first two columns measure containment of the SDP's affine subspace in the hyperplane $S^{\perp}$. The last denotes the minimum eigenvalue of the reducing certificate $S$.} \label{tab:dimacs3} 
         
\end{table}


\section{Conclusion}
We  presented a general technique for  facial reduction that utilizes  approximations of the positive semidefinite cone.  The technique
 is effective on examples arising in practice and for simple  approximation is a  practical pre-processing routine for SDP solvers.  An implementation
 has been made available.  We also gave a post-processing procedure for dual solution recovery that applies generally to cone programs pre-processed using facial reduction. 
 This recovery procedure always succeeds when the cone is polyhedral, but may fail otherwise, illustrating an interesting  difference
between linear programming and optimization over non-polyhedral cones.

\paragraph{Acknowledgments}
The authors thanks Mark Tobenkin for many helpful discussions, and Michael Posa and Anirudha Majumdar for testing early versions of the implemented algorithms.
Michael Posa also provided the original (i.e. unreduced) SDP for Example~\ref{ex:lyap}. The authors thank Cynthia Vinzant for providing the V{\'a}mos matroid  example $\mathcal{V}_{10}$
and  reference \cite{branden2007polynomials}. We thank G\'abor Pataki and Johan L\"ofberg  for  helpful comments.

\bibliographystyle{abbrvnat}
\bibliography{bib}
\end{document}